\documentclass[11pt]{article}
\usepackage{mathrsfs}
\usepackage{amssymb}
\usepackage{amsmath}
\usepackage{amsbsy}
\usepackage{epsfig}
\usepackage{enumerate}
\usepackage{bm}
\usepackage{color}
\newtheorem{lemma}{Lemma}[section]
\newtheorem{theorem}{Theorem}[section]

\newtheorem{corollary}{Corollary}[section]
\newtheorem{definition}{Definiton}[section]
\newtheorem{remark}{Remark}[section]

\newbox\TempBox \newbox\TempBoxA

\newcommand{\red}[1]{\textcolor{red}{#1}}

\def\ep{\textsf{E}} 
\def\Sbep{\widehat{\mathbb E}} 
\def\cSbep{\widehat{\mathcal E}} 
\def\Capc{\mathbb V} 
\def\cCapc{\mathcal V} 

\def\underwiggle 1{
\ifmmode\setbox\TempBox=\hbox{$ 1$}\else\setbox\TempBox=\hbox{
1}\fi \setbox\TempBoxA=\hbox to \wd\TempBox{\hss\char'176\hss}
\rlap{\copy\TempBox}\smash{\lower9pt\hbox{\copy\TempBoxA}} }
\renewcommand{\baselinestretch}{1.7}

\begin{document}

\thispagestyle{empty}
\noindent
The Paper is published at {\em Communications in Mathematics and Statistics} (2015) 3:187-214

 \noindent DOI:10.1007/s40304-015-0055-0

\noindent The red characters are the revisions of the typos.

\begin{center}
 { \LARGE\bf Donsker's invariance principle under the sub-linear expectation with an application to Chung's law of the iterated logarithm$^{\ast}$}
\end{center}

\begin{center} {\sc
Li-Xin Zhang\footnote{Research supported by Grants from the National Natural Science Foundation of China (No. 11225104)  and the Fundamental Research Funds for the Central Universities.
}
}\\
{\sl \small Department of Mathematics, Zhejiang University, Hangzhou 310027} \\
(Email:stazlx@zju.edu.cn)
\end{center}

\renewcommand{\abstractname}{~}
\begin{abstract}
\centerline{\bf Abstract}

We prove a new Donsker's invariance principle for independent and identically distributed random variables under the sub-linear expectation. As applications, the small deviations and Chung's law of the iterated logarithm are obtained.

{\bf Keywords:} sub-linear expectation; capacity;
 central limit theorem; invariance principle; Chung's law of the iterated logarithm; small deviation

{\bf AMS 2010 subject classifications:} 60F15; 60F05; 60H10; 60G48
\end{abstract}

\baselineskip 22pt

\renewcommand{\baselinestretch}{1.7}



\section{ Introduction}\label{sectIntro}
\setcounter{equation}{0}

Let $\{X_n; n\ge 1\}$ be a sequence of independent and identically distributed random variables on a probability space $(\Omega, \mathscr{F}, P)$ such that $EX_1=0$ and $EX_1^2=\sigma^2$. Set $S_n=\sum_{j=1}^n X_j$. In his classical paper Chung (1948) proved the following remarkable result: if $E|X_1|^3<\infty$, then
\begin{equation}\label{eqChung} P\left(\liminf_{n\to \infty} \sqrt{\frac{8\log\log n}{n \pi^2}}\max_{k\le n}|S_k|=  \sigma \right)=1.
\end{equation}
The result (\ref{eqChung}) is  refereed to as the other law of the iterated logarithm or Chung's law of the iterated logarithm, in contrast to the Hartman-Winter law:
\begin{equation}\label{HWLIL} P\left(\limsup_{n\to \infty} \frac{\max_{k\le n}|S_k|}{\sqrt{2n\log\log n}}
=\limsup_{n\to \infty} \frac{|S_n|}{\sqrt{2n\log\log n}}=  \sigma \right)=1.
\end{equation}
 Jian and Pruitt (1975) were the first to prove that (\ref{eqChung}) is still true of only the existence of second moments is assumed. The main tool of Jian and Pruitt (1975) is the following inequality: for any $c<\sigma \pi/\sqrt{8}$, there is a $\eta>1$ such that
\begin{equation}\label{eqJP} P\left( \max_{i\le n}|S_i|< c\Big(\frac{n}{\log\log n}\Big)^{1/2}\right)\le C\frac{1}{(\log n)^{\eta}}
\end{equation}
for all $n$ sufficient large. A more general inequality of this type is the following  small deviation   obtained by Mogul'ski$\check{{\i}}$ (1974):
 if $0<x_n\to 0$ and $n^{1/2}x_n\to \infty$, then
\begin{equation}\label{eqMog} \log P\big(\max_{i\le n}|S_i|\le n^{1/2} x_n\big)\sim -\frac{\pi^2}{8 x_n^2}.
\end{equation}
 A  small deviation theorem and Chung' law of the iterated logarithm for general independent but non-identically distributed random variables was established  by Shao (1995) under a uniform Lindeberg's condition. The key for establishing the small deviations as (\ref{eqJP}),(\ref{eqMog}) and Shao's is the remarkable Donsker's invariance principle as follows:
 let
 $$ W_n(t)=\frac{S_{[nt]}}{\sqrt{n}}+ (nt-[nt])\frac{X_{[nt]+1}}{\sqrt{n}}, \;\; 0\le t\le 1, $$
 then  $W_n\overset{d}\to \sigma W$ in $C[0,1]$ in the sense that
\begin{equation}\label{eqDonsker}  E_P\left[\varphi(W_n)\right]\to E_P\left[\varphi(\sigma W)\right]
\end{equation}
 for any bounded continuous map $\varphi$ from $C[0,1]$ to $\mathbb R$, where $W$ is a standard Browian motion under $P$ and $C[0,1]$ is the space of all continuous function $x:[0,1]\to \mathbb R$ (c.f, Donsker (1951), Billingsley (1968)).

The key in the proof of
the these classical results is the additivity of the probability and the expectation. Under the sub-linear expectation, the Hartman-Winter law of the iterated logarithm  were recently established by Chen and Hu (2014) for bounded random variables, and large deviations were derived by Gao and Xu (2011, 2012).   The main purpose of this paper is to show that (\ref{eqDonsker}) is still true under the sub-linear expectation  and to establish the small deviations similar to (\ref{eqMog}) and Chung's law of the iterated logaritm similar to  (\ref{eqChung}) under the capacities related to the sub-linear expectation.

The general framework of the sub-linear expectation of random variables  in a general function space was introduced by Peng (2006, 2008a, 2008b)
and is a natural extension of the classical linear expectation with the linear property
being replaced by the sub-additivity
and positive homogeneity (c.f  Definition~\ref{def1.1} below).
This simple generalization provides
a  very flexible framework to model non-additive probability and   expectation problems. Take the hedge pricing in fiance  as an example. The famous Black-Shores's formula states that, if a market is
complete and self-financial, then there exists a neutral probability measure $P$ such
that the pricing of any discounted contingent claim $\xi$ in this market is given
by $E_P[\xi]$. However, if the
market is incomplete, such a neutral probability measure is no longer unique, but a
set $\mathscr{P}$ of probability measures. In that case, one can define superhedge pricing $ \Sbep[\xi]=\sup_{Q\in \mathscr{P}} E_Q[\xi]$ and subhedge pricing
$\cSbep[\xi]=-\Sbep[-\xi]=\inf_{Q\in \mathscr{P}} E_Q[\xi]$, respectively. Then $\Sbep[\xi]$ is sub-linear,   as a functional operator of random
variables, and the related  capacity $\Capc(A)=\sup_{Q\in \mathscr{P}} Q(A)$ is non-additive, as a function operator of events. The extension of linear expectations to sub-linear expectations also  produces many interesting properties different
 from the classic ones.
  For examples,  the limit in the law of large numbers is no longer a contact,
  and, comparing to the  classical one-dimensional  normal  distribution which is characterized by the Stein equation, an
ordinary differential equation (ODE),
a normal distribution under the sub-linear expectation is characterized by a time-space parabolic partial  different equation (PDE). Recently, Hu and Li (2014) showed that the characteristic function cannot determine the distribution of
random variables on the sub-linear expectation space.
  Roughly speaking, a sub-linear expectation is related to a group of unknown linear expectations
and the distribution under a sub-linear expectation is related to a group of probabilities (c.f. Lemma 2.4 of Peng (2008b)).
For more properties of the sub-linear expectations, one can refer to Peng (2008b,2009),
where the notion of independent and identically distributed random variables under  the sub-linear expectations was introduced  and a new central limit theorem was estiblished. It is a natural and interesting problem whether (1.4) is true when the linear expectation $E$ is replaced by the sub-linear expectation.

In the classical probability space, (1.4) is equivalent  to the weak convergence of related probability measures in the metric space $C[0,1]$ equipped with the super-metric $\|x-y\|=\sup\limits_{0\le t\le 1}|x(t)-y(t)|$.
Classically,   the weak convergence of probability measures in $C[0,1]$ is showed by verifying the convergence of the finite-dimensional distributions and  the tightness of the probability measures. We will find that this way is also valid for proving   Donsker's invariance principle in the sub-linear expectation space, though 
there is no longer any one-to-one relationship between the convergence of sub-linear expectations and the convergence of related capacities. In the next section, we state basic settings in a sub-linear expectation space including, capacity, independence, identical distribution, G-Brownian motion etc.    The main result  on   Donsker's invariance principle is established in Section \ref{sectMain} by assuming that the convergence of finite-dimensional distributions and the tightness are proved.  And also, Chung's law of the iterated logarithm is established by assuming the small deviations. We consider the convergence of the finite-dimensional distributions in Section \ref{sectFinite}, and consider the tightness in Section \ref{sectTight}. In Section \ref{sectDeviation}, we establish the small deviations by applying   Donsker's invariance principle. All of the results are established of only the existence of second moments is assumed.

\section{Basic Settings}\label{sectBasic}
\setcounter{equation}{0}

We use the framework and notations of Peng (2008b). Let  $(\Omega,\mathcal F)$
 be a given measurable space  and let $\mathscr{H}$ be a linear space of real functions
defined on $(\Omega,\mathcal F)$ such that if $X_1,\ldots, X_n \in \mathscr{H}$  then $\varphi(X_1,\ldots,X_n)\in \mathscr{H}$ for each
$\varphi\in C_b(\mathbb R_n)\bigcup  C_{l,Lip}(\mathbb R_n)$,  where $C_b(\mathbb R_n)$ denote  the
space  of all bounded  continuous functions and $C_{l,Lip}(\mathbb R_n)$ denotes the linear space of (local Lipschitz)
functions $\varphi$ satisfying
\begin{eqnarray*} & |\varphi(\bm x) - \varphi(\bm y)| \le  C(1 + |\bm x|^m + |\bm y|^m)|\bm x- \bm y|, \;\; \forall \bm x, \bm y \in \mathbb R_n,&\\
& \text {for some }  C > 0, m \in \mathbb  N \text{ depending on } \varphi. &
\end{eqnarray*}
$\mathscr{H}$ is considered as a space of ``random variables''. In this case we denote $X\in \mathscr{H}$.
Further, we let $C_{b,Lip}(\mathbb R_n)$  denote  the
space  of all bounded and Lipschitz functions   on $\mathbb R_n$.

\begin{remark} It is easily seen that if $\varphi_1,\varphi_2\in C_{l,Lip}(\mathbb R_n)$, then $\varphi_1\vee \varphi_2, \varphi_1\wedge \varphi_2\in C_{l,Lip}(\mathbb R_n)$ because $\varphi_1\vee \varphi_2=\frac{1}{2}( \varphi_1+ \varphi_2+|\varphi_1- \varphi_2|)$, $\varphi_1\wedge \varphi_2=\frac{1}{2}( \varphi_1+ \varphi_2-|\varphi_1- \varphi_2|)$.
\end{remark}

\subsection{Sub-linear expectation and capacity}
\begin{definition}\label{def1.1} A {\bf sub-linear expectation} $\Sbep$ on $\mathscr{H}$  is a functional $\Sbep: \mathscr{H}\to \overline{\mathbb R}$ satisfying the following properties: for all $X, Y \in \mathscr H$, we have
\begin{description}
  \item[\rm (a)] {\bf Monotonicity}: If $X \ge  Y$ then $\Sbep [X]\ge \Sbep [Y]$;
\item[\rm (b)] {\bf Constant preserving} : $\Sbep [c] = c$;
\item[\rm (c)] {\bf Sub-additivity}: $\Sbep[X+Y]\le \Sbep [X] +\Sbep [Y ]$ whenever $\Sbep [X] +\Sbep [Y ]$ is not of the form $+\infty-\infty$ or $-\infty+\infty$;
\item[\rm (d)] {\bf Positive homogeneity}: $\Sbep [\lambda X] = \lambda \Sbep  [X]$, $\lambda\ge 0$.
 \end{description}
 Here $\overline{\mathbb R}=[-\infty, \infty]$. The triple $(\Omega, \mathscr{H}, \Sbep)$ is called a sub-linear expectation space. Give a sub-linear expectation $\Sbep $, let us denote the conjugate expectation $\cSbep$of $\Sbep$ by
$$ \cSbep[X]:=-\Sbep[-X], \;\; \forall X\in \mathscr{H}. $$
\end{definition}

From the definition, it is easily shown that    $\cSbep[X]\le \Sbep[X]$, $\Sbep[X+c]= \Sbep[X]+c$ and $\Sbep[X-Y]\ge \Sbep[X]-\Sbep[Y]$ for all
$X, Y\in \mathscr{H}$ with $\Sbep[Y]$ being finite. Further, if $\Sbep[|X|]$ is finite, then $\cSbep[X]$ and $\Sbep[X]$ are both finite.

\bigskip

Next, we introduce the capacities corresponding to the sub-linear expectations.
Let $\mathcal G\subset\mathcal F$. A function $V:\mathcal G\to [0,1]$ is called a capacity if
$$ V(\emptyset)=0, \;V(\Omega)=1 \; \text{ and } V(A)\le V(B)\;\; \forall\; A\subset B, \; A,B\in \mathcal G. $$
It is called to be sub-additive if $V(A\bigcup B)\le V(A)+V(B)$ for all $A,B\in \mathcal G$  with $A\bigcup B\in \mathcal G$.

 Let $(\Omega, \mathscr{H}, \Sbep)$ be a sub-linear space,
and  $\cSbep $  be  the conjugate expectation of $\Sbep$. It is natural to define the capacity of a set $A$ to be the sub-linear expectation of
the   indicator function $I_A$ of $A$. However, $I_A$ may be not in $\mathscr{H}$. So, we denote a pair $(\Capc,\cCapc)$ of capacities by
$$ \Capc(A):=\inf\{\Sbep[\xi]: I_A\le \xi, \xi\in\mathscr{H}\}, \;\; \cCapc(A):= 1-\Capc(A^c),\;\; \forall A\in \mathcal F, $$
where $A^c$  is the complement set of $A$.
Then
\begin{equation}\label{eq1.3} \begin{matrix}
&\Capc(A)=\Sbep[I_A], \;\; \cCapc(A)= \cSbep[I_A],\;\; \text{ if } I_A\in \mathscr H\\
&\Sbep[f]\le \Capc(A)\le \Sbep[g], \;\;\cSbep[f]\le \cCapc(A) \le \cSbep[g],\;\;
\text{ if } f\le I_A\le g, f,g \in \mathscr{H}.
\end{matrix}
\end{equation}
It is obvious that $\Capc$ is sub-additive. But $\cCapc$ and $\cSbep$ may  be not sub-additive. However, we have
\begin{equation}\label{eq1.4}
  \cCapc(A\bigcup B)\le \cCapc(A)+\Capc(B) \;\;\text{ and }\;\; \cSbep[X+Y]\le \cSbep[X]+\Sbep[Y]
\end{equation}
due to the fact that $\Capc(A^c\bigcap B^c)=\Capc(A^c\backslash B)\ge \Capc(A^c)-\Capc(B)$ and $\Sbep[-X-Y]\ge \Sbep[-X]-\Sbep[Y]$.

Further, we define an extension of $\Sbep^{\ast}$ of $\Sbep$ by
$$ \Sbep^{\ast}[X]=\inf\{\Sbep[Y]:X\le Y, \; Y\in \mathscr{H}\}, \;\; \forall X:\Omega\to \mathbb R, $$
where $\inf\emptyset=+\infty$. Then
$$ \begin{matrix}
&\Sbep^{\ast}[X]=\Sbep[X]\; \text{ if } X\in \mathscr H, \;\; \;\; \Capc(A)=\Sbep^{\ast}[I_A], \\
&\Sbep[f]\le \Sbep^{\ast}[X]\le \Sbep[g]\;\;
\text{ if } f\le X\le g, f,g \in \mathscr{H}.
\end{matrix}
$$

 \begin{definition}\label{def3.1}
\begin{description}
\item{\rm (I)} A sub-linear expectation $\Sbep: \mathscr{H}\to \mathbb R$ is called to be  countably sub-additive if it satisfies
\begin{description}
  \item[\rm (e)] {\bf Countable sub-additivity}: $\Sbep[X]\le \sum_{n=1}^{\infty} \Sbep [X_n]$, whenever $X\le \sum_{n=1}^{\infty}X_n$,
  $X, X_n\in \mathscr{H}$ and
  $X\ge 0, X_n\ge 0$, $n=1,2,\ldots$;
 \end{description}
It is called to continuous if it satisfies
\begin{description}
  \item[\rm (f) ]  {\bf Continuity from below}: $\Sbep[X_n]\uparrow \Sbep[X]$ if $0\le X_n\uparrow X$, where $X_n, X\in \mathscr{H}$;
  \item[\rm (g) ] {\bf Continuity from above}: $\Sbep[X_n]\downarrow \Sbep[X]$ if $0\le X_n\downarrow X$, where $X_n, X\in \mathscr{H}$.
\end{description}
\item{\rm (II)}  A function $V:\mathcal F\to [0,1]$ is called to be  countably sub-additive if
$$ V\Big(\bigcup_{n=1}^{\infty} A_n\Big)\le \sum_{n=1}^{\infty}V(A_n) \;\; \forall A_n\in \mathcal F. $$

\item{\rm (III)}  A capacity $V:\mathcal F\to [0,1]$ is called a continuous capacity if it satisfies
\begin{description}
  \item[\rm (III1) ] {\bf Ccontinuity from below}: $V(A_n)\uparrow V(A)$ if $A_n\uparrow A$, where $A_n, A\in \mathcal F$;
  \item[\rm (III2) ] {\bf Continuity from above}: $V(A_n)\downarrow  V(A)$ if $A_n\downarrow A$, where $A_n, A\in \mathcal F$.
\end{description}
\end{description}
\end{definition}

It is obvious that a continuous sub-additive capacity $V$ (resp. a sub-linear expectation $\Sbep$) is countably sub-additive.

\subsection{Independence and distribution}
\begin{definition} ({\em Peng (2006, 2008b)})

\begin{description}
  \item[ \rm (i)] ({\bf Identical distribution}) Let $\bm X_1$ and $\bm X_2$ be two $n$-dimensional random vectors defined
respectively in sub-linear expectation spaces $(\Omega_1, \mathscr{H}_1, \Sbep_1)$
  and $(\Omega_2, \mathscr{H}_2, \Sbep_2)$. They are called identically distributed, denoted by $\bm X_1\overset{d}= \bm X_2$  if
$$ \Sbep_1[\varphi(\bm X_1)]=\Sbep_2[\varphi(\bm X_2)], \;\; \forall \varphi\in C_{l,Lip}(\mathbb R_n), $$
whenever the sub-expectations are finite. A sequence $\{X_n;n\ge 1\}$ of random variables is said to be identically distributed if $X_i\overset{d}= X_1$ for each $i\ge 1$.
\item[\rm (ii)] ({\bf Independence})   In a sub-linear expectation space  $(\Omega, \mathscr{H}, \Sbep)$, a random vector $\bm Y =
(Y_1, \ldots, Y_n)$, $Y_i \in \mathscr{H}$ is said to be independent to another random vector $\bm X =
(X_1, \ldots, X_m)$ , $X_i \in \mathscr{H}$ under $\Sbep$  if for each test function $\varphi\in C_{l,Lip}(\mathbb R_m \times \mathbb R_n)$
we have
$$ \Sbep [\varphi(\bm X, \bm Y )] = \Sbep \big[\Sbep[\varphi(\bm x, \bm Y )]\big|_{\bm x=\bm X}\big],$$
whenever $\overline{\varphi}(\bm x):=\Sbep\left[|\varphi(\bm x, \bm Y )|\right]<\infty$ for all $\bm x$ and
 $\Sbep\left[|\overline{\varphi}(\bm X)|\right]<\infty$.
 \item[\rm (iii)] ({\bf IID random variables}) A sequence of random variables $\{X_n; n\ge 1\}$
 is said to be independent and identically distributed (IID), if
 $X_i\overset{d}=X_1$ and $X_{i+1}$ is independent to $(X_1,\ldots, X_i)$ for each $i\ge 1$.
 \end{description}
\end{definition}

\subsection{G-normal distribution and G-Brownian motion}
Let $0\le \underline{\sigma}\le \overline{\sigma}<\infty$ and $G(\alpha)=\frac{1}{2}(\overline{\sigma}^2 \alpha^+ - \underline{\sigma}^2 \alpha^-)$. $X$ is call a normal $N\big(0, [\underline{\sigma}^2, \overline{\sigma}^2]\big)$ distributed random variable  (write $X\sim N\big(0, [\underline{\sigma}^2, \overline{\sigma}^2]\big)$) under $\Sbep$, if for any $\varphi\in C_b(\mathbb R)$, the function $u(x,t)=\Sbep\left[\varphi\left(x+\sqrt{t} X\right)\right]$ ($x\in \mathbb R, t\ge 0$) satisfies the following heat equation:
      $$ \partial_t u -G\left( \partial_{xx}^2 u\right) =0, \;\; u(0,x)=\varphi(x). $$

 Let $C[0,1]$ be a function space of continuous functions   on $[0,1]$ equipped with the super-norm $\|x\|=\sup\limits_{0\le t\le 1}|x(t)|$ and $C_b\big(C[0,1]\big)$ is the set of bounded continuous  functions $h(x):C[0,1]\to \mathbb R$.The modulus of the continuity of an element   $x\in C[0,1]$  is defined by $$\omega_{\delta}(x)=\sup_{|t-s|<\delta}|x(t)-x(s)|.$$

 Denis, Hu and Peng (2011) showed that there is a sub-linear expectation space $\big(\widetilde{\Omega}, \widetilde{\mathscr{H}},\widetilde{\mathbb E}\big)$ with
$\widetilde{\Omega}= C[0,1]$ and $C_b\big(C[0,1]\big)\subset \widetilde{\mathscr{H}}$ such that $\widetilde{\mathbb E}$ is countably sub-additive, and
the canonical process $W(t)(\omega) = \omega_t  (\omega\in \widetilde{\Omega})$ is a G-Brownian motion with
$W(1)\sim N\big(0, [\underline{\sigma}^2, \overline{\sigma}^2]\big)$ under $\widetilde{\mathbb E}$, i.e.,
  for all $0\le t_1<\ldots<t_n\le 1$, $\varphi\in C_{b,lip}(\mathbb R^n)$,
\begin{equation}\label{eqBM} \widetilde{\mathbb E}\left[\varphi\big(W(t_1),\ldots, W(t_{n-1}), W(t_n)-W(t_{n-1})\big)\right]
  =\widetilde{\mathbb E}\left[\psi\big(W(t_1),\ldots, W(t_{n-1})\big)\right],
  \end{equation}
  where $\psi(x_1,\ldots, x_{n-1})\big)=\widetilde{\mathbb E}\left[\varphi\big(x_1,\ldots, x_{n-1}, \sqrt{t_n-t_{n-1}}W(1)\big)\right]$.

Denis, Hu and Peng (2011) also showed the following representation of the G-Brownian motion (c.f, Theorem 52).
 \begin{lemma} \label{DenisHuPeng}
 Let $(\Omega, \mathscr{F}, P) $ be a probability
measure space and $\{B(t)\}_{t\ge 0}$  is a $P$-Brownian motion. Then for all $\varphi\in C_b(\widetilde{\Omega})$,
$$ \widetilde{\mathbb E}\left[\varphi(W)\right]=\sup_{\theta\in \Theta}\ep_P\left[\varphi\left(W_{\theta}\right)\right],\;\;
W_{\theta}(t) = \int_0^t\theta(s) dB(s), $$
where
\begin{eqnarray*}
&\Theta=\left\{ \theta:\theta(t) \text{ is } \mathscr{F}_t\text{-adapted process such that }  \underline{\sigma}\le \theta(t)\le \overline{\sigma}\right\},&\\
& \mathscr{F}_t=\sigma\{B(s):0\le s\le t\}\vee \mathscr{N}, \;\; \mathscr{N} \text{ is the collection of } P\text{-null subsets}. &
\end{eqnarray*}
\end{lemma}

We   we denote a pair of   capacities corresponding to the sub-linear expectation $\widetilde{\mathbb E}$ by  $(\widetilde{\Capc},\widetilde{\cCapc})$, and the extension of $\widetilde{\mathbb E}$ by  $\widetilde{\mathbb E}^{\ast}$.   By using Lemma \ref{DenisHuPeng}, one can show that
$$ \widetilde{\mathbb E}^{\ast}\left[\sup_{0\le t\le 1}|W(t)|^p\right]=\overline{\sigma}^p
E_P\left[\sup_{0\le t\le 1}|B(t)|^p\right] $$
and for $x\ge 0$
\begin{equation}\label{eqBrownian1}
  \widetilde{\Capc}\left(\sup_{0\le t\le 1}|W(t)|\ge x\right)=   P\left(\sup_{0\le t\le 1}\overline{\sigma}|B(t)|\ge x\right),
  \end{equation}
 $$ \widetilde{\cCapc}\left(\sup_{0\le t\le 1}|W(t)|\ge x\right)=   P\left(\sup_{0\le t\le 1}\underline{\sigma}|B(t)|\ge x\right),   $$
$$ \widetilde{\Capc}\left(\sup_{0\le t\le 1}W(t)\ge x\right)=   2P\left( \overline{\sigma} B(1) \ge x\right), $$
$$\widetilde{\cCapc}\left(\sup_{0\le t\le 1} W(t) \ge x\right)=   2P\left(   \underline{\sigma} B(1)\ge x\right)$$
 (c.f, Lemma \ref{lemma7}  below).

 In the sequel of this paper, the sequence $\{X_n;n\ge 1\}$ of the random variables are considered in $(\Omega, \mathscr{H}, \Sbep)$ and Brownian motions are considered in $(\widetilde{\Omega}, \widetilde{\mathscr{H}}, \widetilde{\mathbb E})$. We suppose $\{X_n; n\ge 1\}$ is a sequence of independent and identically distributed random variables in $(\Omega, \mathscr{H}, \Sbep)$  with  $\Sbep[X_1]=\cSbep[X_1]=0$, $\Sbep[X_1^2]=\overline{\sigma}^2$ and $\cSbep[X_1^2]=\underline{\sigma}^2$, and suppose    $W(t)$ is a G-Brownian motion on $(\widetilde{\Omega}, \widetilde{\mathscr H}, \widetilde{\mathbb E})$  with $W(1)\sim N(0,[\underline{\sigma}^2,\overline{\sigma}^2])$. Denote $S_0=0$, $S_n=\sum_{k=1}^n X_k$.

\section{Main results}\label{sectMain}
\setcounter{equation}{0}

 Define the $C[0,1]$-valued random variable $W_n$ by setting
$$W_n(t)=\begin{cases}
 S_k/\sqrt{n},  \; \text{ if }  t=k/n \; (k=0,1,\ldots, n);\\
 \text{ extended by linear interpolation in each interval }\\
  \qquad \quad \big[[k-1]n^{-1}, kn^{-1}\big].
\end{cases}$$

Our  first result is   the following Donsker's invariance principle, or called  the functional central limit theorem.
\begin{theorem}\label{thFCLT} Suppose $\Sbep[(X_1^2-b)^+]\to 0$ as $b\to \infty$. Then for all bounded   continuous function $\varphi:C[0,1]\to \mathbb R$,
\begin{equation}\label{eqthFCLT.1}
\Sbep\left[\varphi(W_n)\right]\to \widetilde{\mathbb E}\left[\varphi(W)\right].
 \end{equation}
\end{theorem}

 \begin{remark} Note $C_b\left(C[0,1]\right)\subset\widetilde{\mathscr{H}}$. The expectation on the right hands of (\ref{eqthFCLT.1}) is well-defined. On the other hand, the map $g: (S_1/\sqrt{n},\ldots,S_n/\sqrt{n})\to W_n(t)$ is a  continuous one-to-one map of $(S_1/\sqrt{n},\ldots,S_n/\sqrt{n})$. So, $\varphi\circ g\in C_b(\mathbb R^n)$. It follows that $\varphi(W_n)\in \mathscr{H}$. So, even though $\Sbep$ has no definition on $C_b\left(C[0,1]\right)$, the expectation  on the left hand of (\ref{eqthFCLT.1}) is well-defined.
  \end{remark}

  There are there immediate corollaries of Theorem \ref{thFCLT}.
   \begin{corollary}\label{cor1} Suppose $\Sbep[(X_1^2-b)^+]\to 0$ as $b\to \infty$. If $h_n(x,y)$ and $h(x,y)$: $C[0,1]\times C[0,1]\to \mathbb R$ are bounded continuous functions for which
 \begin{eqnarray*}
 & |h_n(x,y)|\le M, \;\; |h(x,y)|\le M,& \\
 & h_n(x_n,y_n)\to h(x,y) \text{ whenever }  x_n\to x, y_n\to y,
 \end{eqnarray*}
 then
 \begin{equation}\label{eqthFCLT.2}
\Sbep\left[h_n(W_n, y)\right]\to  \widetilde{\mathbb E}\left[h(W,y)\right]
\;\; \text{ uniformly  in } y \in K,
 \end{equation}
 for any compact set $K\subset C[0,1]$.
\end{corollary}

\begin{corollary}\label{cor2} Suppose $p\ge 2$ and $\Sbep[(|X_1|^p-b)^+]\to 0$ as $b\to \infty$. Then for all continuous function $\varphi:C[0,1]\to \mathbb R$ with $|\varphi(x)|\le C(1+\|x\|^p)$,
\begin{equation}\label{eqcor2.1}
\Sbep^{\ast}\left[\varphi(W_n)\right]\to \widetilde{\mathbb E}^{\ast}\left[\varphi(W)\right]
 \end{equation}
 where $\Sbep^{\ast}$ and $\widetilde{\mathbb E}^{\ast}$ are extensions of $\Sbep$ and $\widetilde{\mathbb E}$, respectively. In particular,
 $$ \Sbep\left[\max_{k\le n}\left|\frac{S_k}{\sqrt{n}}\right|^p\right]\to \overline{\sigma}^p E_P[\sup_{0\le t\le 1}|B(t)|^p],\;\; \cSbep\left[\max_{k\le n}\left|\frac{S_k}{\sqrt{n}}\right|^p\right]\to \underline{\sigma}^p E_P[\sup_{0\le t\le 1}|B(t)|^p]. $$
    \end{corollary}

  \begin{corollary}\label{cor3} Suppose $\Sbep[(X_1^2-b)^+]\to 0$ as $b\to \infty$. Then for all $x>0$,
  $$\Capc\left(\max_{k\le n}|S_k|/\sqrt{n}\ge x\right)\to P\left(\sup_{0\le t\le 1}\overline{\sigma}|B(t)|\ge x\right), $$
$$\cCapc\left(\max_{k\le n}|S_k|/\sqrt{n}|\ge x\right)\to P\left(\sup_{0\le t\le 1}\underline{\sigma}|B(t)|\ge x\right), $$
$$\Capc\left(\max_{k\le n} S_k /\sqrt{n}\ge x\right)\to 2P\left( \overline{\sigma} B(1) \ge x\right), $$
$$\cCapc\left(\max_{k\le n} S_k /\sqrt{n} \ge x\right)\to 2P\left(   \underline{\sigma} B(\red{1}) \ge x\right). $$
  \end{corollary}

Note
 $$ \lim_{x\to 0+} x^2\log P\left(\sup_{0\le t\le 1}|B(t)|\le x\right)
 =-\frac{\pi^2}{8}. $$
 From Corollary \ref{cor3}, we conclude that
 $$ \lim_{x\to 0+} x^2\lim_{n\to \infty} \log \Capc\left(\max_{k\le n}|S_k|/\sqrt{n}\le x\right)=-\frac{\pi^2\underline{\sigma}^2}{8}, $$
 $$ \lim_{x\to 0+} x^2\lim_{n\to \infty}\log \cCapc\left(\max_{k\le n}|S_k|/\sqrt{n}\le x\right)=-\frac{\pi^2\overline{\sigma}^2}{8}. $$
In Section \ref{sectDeviation}, we will prove  more accurate results refereed to as the small deviations:
\begin{align}\label{eqDeviation0.1}
& \log \Capc\left(\max_{k\le n}|S_k|\le n^{1/2}x_n\right)\sim -\frac{\pi^2 \underline{\sigma}^2}{8x_n^2},\\
\label{eqDeviation0.2}
& \log \cCapc\left(\max_{k\le n}|S_k|\le n^{1/2}x_n\right)\sim -\frac{\pi^2 \overline{\sigma}^2}{8x_n^2}
\end{align}
whenever $0<x_n\to 0$ and $\sqrt{n}x_n\to \infty$. By the small deviations, we obtain  our another main result which gives  Chung's law of the iterated logarithm.
 \begin{theorem}\label{thChungLIL} Suppose $\Sbep[(X_1^2-b)^+]\to 0$ as $b\to \infty$ and $\Capc$ is continuous. Then
 \begin{equation}\label{eqChungLIL.1}
 \cCapc\left(\underline{\sigma}\le \liminf_{n\to \infty} \sqrt{\frac{8\log\log n}{n \pi^2}}\max_{k\le n}|S_k|\le \overline{\sigma}\right)=1
 \end{equation}
 and
 \begin{equation}\label{eqChungLIL.2}
 \Capc\left(  \liminf_{n\to \infty} \sqrt{\frac{8\log\log n}{n \pi^2}}\max_{k\le n}|S_k|= \underline{\sigma}\right)=1.
 \end{equation}
 \end{theorem}

 \bigskip
Next, we give a sketch of the proof of Theorems \ref{thFCLT}, \ref{thChungLIL} and Corollaries \ref{cor1}- \ref{cor3}. We need  the following Rosenthal type inequalities under $\Sbep$ which have been obtained  by Zhang (2016).
 \begin{lemma}\label{thRIneq} ({\em Rosnethal's inequality})
   Let $\{X_1,\ldots, X_n\}$ be a sequence  of  independent random variables
 in $(\Omega, \mathscr{H}, \Sbep)$.
    Then
\begin{align}\label{eqthRIneq.2}
\Sbep\left[\max_{k\le n} \left|S_k\right|^p\right]\le  & C_p\left\{ \sum_{k=1}^n \Sbep [|X_k|^p]+\left(\sum_{k=1}^n \Sbep [|X_k|^2]\right)^{p/2} \right. \nonumber
 \\
& \qquad \left. +\left(\sum_{k=1}^n \big\{\big(\cSbep [X_k]\big)^-+\big(\Sbep [X_k]\big)^+\big\}\right)^{p}\right\}, \;\; \text{ for }   p\ge 2.
\end{align}
In particular, if $\Sbep[X_k]=\Sbep[-X_k]= 0$, $k=1,\ldots, n$, then
$$
\Sbep\left[\max_{k\le n}|S_k|^p \right]\le C_p\left\{ \sum_{k=1}^n \Sbep [|X_k|^p]+\left(\sum_{k=1}^n \Sbep [|X_k|^2]\right)^{p/2}\right\},
\;\; \text{ for }   p\ge 2.
$$
\end{lemma}

\bigskip

{\bf Proof of Theorem \ref{thFCLT}.} The proof is based on the convergence of the finite-dimensional distributions of $W_n$  under $\Sbep$ and the tightness of $W_n$ under $\Capc$  which are given in Section \ref{sectFinite} and Section \ref{sectTight}, respectively. Here, we give the proof of theorem after assuming these results.

 For $0< t_1<t_2\ldots<t_d\le 1$, we define the projection $\pi_{t_1,\ldots,t_d}$ from $C[0,1]$ to $\mathbb R^d$ by
$$ \pi_{t_1,\ldots,t_d}x=(x(t_1),\ldots, x(t_d)), $$
and define a map $\Pi^{-1}_{t_1,\ldots,t_d}$ from $\mathbb R^d$ to $C[0,1]$  by
$$\Pi^{-1}_{t_1,\ldots,t_d}(x_1,\ldots,x_d)= \begin{cases}
0,  \; \text{ if }  t=0;\;\;
 x_k,  \; \text{ if }  t=t_k \; (k=0,1,\ldots, n);\\
 \text{ extended by linear interpolation in each interval }\\
  \qquad \quad \big[t_{k-1}, t_k\big].
\end{cases}$$
Then $\pi_{t_1,\ldots,t_d}$ and $\Pi_{t_1,\ldots,t_d}^{-1}$ are both continuous maps. Denote $\widetilde{\pi}_{t_1,\ldots,t_d}=\Pi^{-1}_{t_1,\ldots,t_d}\circ \pi_{t_1,\ldots,t_d}$.
Then $\widetilde{\pi}_{t_1,\ldots,t_d}:C[0,1]\to C[0,1]$ is continuous and
$$ W_n=\Pi_{1/n, 2/n,\ldots,n/n}^{-1}\big(S_1/\sqrt{n},S_2/\sqrt{n},\ldots, S_n/\sqrt{n}\big). $$
Let $\varphi\in C_b\big(C[0,1])$. Then $\varphi(\widetilde{\pi}_{t_1,\ldots,t_d} x)=\varphi\circ\Pi^{-1}_{t_1,\ldots,t_d}(x(t_1),\ldots, x(t_d))$ and $\varphi\circ\Pi^{-1}_{t_1,\ldots,t_d}\in C_b(\mathbb R^d)$. By Theorem \ref{thfinite} on the convergence of the finite-dimensional distributions of $W_n$, it follows that
\begin{align*}
&\lim_{n\to \infty}\Sbep\left[\varphi\left(\widetilde{\pi}_{t_1,\ldots,t_d} W_n\right)\right]
=\lim_{n\to \infty} \Sbep\left[\varphi\circ\Pi^{-1}_{t_1,\ldots,t_d}\left(W_n(t_1),\ldots,W_n(t_d)\right)\right]\\
& \qquad =\Sbep\left[\varphi\circ\Pi^{-1}_{t_1,\ldots,t_d}\left(W(t_1),\ldots,W(t_d)\right)\right]
=\Sbep\left[\varphi\left(\widetilde{\pi}_{t_1,\ldots,t_d} W\right)\right].
\end{align*}
Now, let $t_0=0$, $t_{d+1}=1$, and suppose that $t_{i+1}-t_i<\delta$ for $i=0,\ldots, d$. Recall  $\omega_{\delta}(x)=\sup\limits_{|t-s|<\delta}|x(t)-x(s)|$ and $\|x\|=\sup\limits_{0\le t\le 1}|x(t)|$. It is easily seen that
$  \left\|\widetilde{\pi}_{t_1,\ldots,t_d}x -x \right\|\le \omega_{\delta}(x). $
Let $\epsilon>0$ be given. Since  $\varphi$ is a continuous function, for each $x$, there is an $\epsilon_x>0$ such that
$$ \left|\varphi(x)-\varphi(y)\right|<\epsilon \text{ whenever } \|y-x\|<\epsilon_x. $$
Let $K\subset C[0,1]$ is a compact set. Then it can be covered by a union of finite many of the sets $\{y: \|y-x\|<\epsilon_x\}$, $x\in K$. So, there is an $\epsilon_K>0$ such that
 $ \left|\varphi(x)-\varphi(y)\right|<\epsilon$  whenever $\|y-x\|<\epsilon_K$ and $x\in K$.
Denote $M=\sup_x|\varphi(x)|$. It follows that
 \begin{align*}
 \left|\varphi\left(\widetilde{\pi}_{t_1,\ldots,t_d} x\right)-\varphi(x)\right|
 <\epsilon+2M I\{\omega_{\delta}(x)\ge \epsilon_K\}+2MI\{x\not\in K\}.
 \end{align*}
  By Theorems \ref{thTight} and \ref{thTight2} on the tightness of $\{W_n\}$ and $W$, respectively,  we can choose $K$ and $\delta$ such that
 $$ \sup_n\Capc\left(\omega_{\delta}(W_n)\ge \epsilon_K\right)+  \sup_n\Capc\left(W_n\not\in K\right)
 \le \frac{\epsilon}{4M} \;\; \text{and }$$
 $$\widetilde{\Capc}\left(\omega_{\delta}(W)\ge \epsilon_K\right)+  \widetilde{\Capc}\left(W\not\in K\right)
 \le \frac{\epsilon}{4M}. $$
 Hence
 \begin{align*}
 & \left| \Sbep\left[ \varphi(W_n)\right]- \widetilde{\mathbb E}\left[ \varphi(W)\right]\right|\\
 \le & \left| \Sbep\left[ \varphi\big(\widetilde{\pi}_{t_1,\ldots, t_d}W_n\big)\right]- \widetilde{\mathbb E}\left[\varphi\big(\widetilde{\pi}_{t_1,\ldots, t_d}W\big)\right]\right| \\
 &+\left|\Sbep\left[\varphi(W_n)\right]-\Sbep\left[\varphi\left(\widetilde{\pi}_{t_1,\ldots,t_d} W_n\right)\right]\right|
 +\left|\widetilde{\mathbb E}\left[\varphi\big(\widetilde{\pi}_{t_1,\ldots, t_d}W\big)\right]-\widetilde{\mathbb E}\left[ \varphi(W)\right]\right|\\
 \le & \left| \Sbep\left[ \varphi\big(\widetilde{\pi}_{t_1,\ldots, t_d}W_n\big)\right]- \widetilde{\mathbb E}\left[\varphi\big(\widetilde{\pi}_{t_1,\ldots, t_d}W\big)\right]\right| \\
 &+2\epsilon+ 2M \Capc\left(\omega_{\delta}(W_n)\ge \epsilon_K\right)+2M \Capc\left(W_n\not\in K\right) \\
 &+ 2M \widetilde{\Capc}\left(\omega_{\delta}(W)\ge \epsilon_K\right)+2M \widetilde{\Capc}\left(W\not\in K\right)\\
 \le & \left| \Sbep\left[ \varphi\big(\widetilde{\pi}_{t_1,\ldots, t_d}W_n\big)\right]- \widetilde{\mathbb E}\left[\varphi\big(\widetilde{\pi}_{t_1,\ldots, t_d}W\big)\right]\right|+3\epsilon.
 \end{align*}
Letting $n\to \infty$ and then $\epsilon\to 0$ completes the proof of (\ref{eqthFCLT.1}). $\Box$

\begin{remark} Here we give a direct proof of Theorem \ref{thFCLT}. In Section \ref{sectTight}, we  will give another proof by using the results of Peng (2010).

\end{remark}

\bigskip
{\bf Proof of Corollary \ref{cor1}.}
Let $\{\delta_n\}$ be a sequence of positive numbers with $\delta_n\downarrow 0$, and
let $y_n\in K$ such that
$$ \left|\Sbep\left[h_n(W_n, y_n)\right]-  \widetilde{\mathbb E}\left[h(W,y_n)\right]\right|
\ge \sup_{y \in K}\left| \Sbep\left[h_n(W_n, y)\right]- \widetilde{\mathbb E}\left[h(W,y)\right]\right|-\delta_n.
 $$
Since $K$ is compact, any subsequence of $\{y_n\}$ has a further convergent subsequence. Without loss of generality, we assume $y_n\to y\in K$. Let $g_n(x)=h_n(x,y_n)$ and $g(x)=h(x,y)$. Then
$ g_n(x_n)\to g(x) $   whenever $ x_n \to x,$
which implies that for any compact set $K_1\subset C[0,1]$,
$\sup_{x\in K_1}|g_n(x)-g(x)|\to 0. $
For any $\epsilon>0$, by Theorem  \ref{thTight}, one can choose $K_1$ such that
$ \sup_n \Capc\left(W_n\not\in K_1\right) <\epsilon/(2M). $
It follows that
\begin{align*}
&\sup_{y \in K}\left| \Sbep\left[h_n(W_n, y)\right]-  \widetilde{\mathbb E}\left[h(W,y)\right]\right|
\le \left| \Sbep\left[g_n(W_n)\right]- \widetilde{\mathbb E}\left[g(W)\right]\right|+\delta_n \\
\le & \left| \Sbep\left[g_n(W_n)\right]- \Sbep\left[g(W_n)\right]\right|
+\left| \Sbep\left[g(W_n)\right]- \widetilde{\mathbb E}\left[g(W)\right]\right|+\delta_n\\
\le & \left| \Sbep\left[g(W_n)\right]-  \widetilde{\mathbb E}\left[g(W)\right]\right|
+\sup_{x\in K_1}|g_n(x)-g(x)|+2M\cdot \epsilon/(2M)+\delta_n.
\end{align*}
From  (\ref{eqTight.1}), it follows that
$$\limsup_{n\to \infty}\sup_{y \in K}\left| \Sbep\left[h_n(W_n, y)\right]-  \widetilde{\mathbb E}\left[h(W,y)\right]\right|\le \epsilon. $$
The proof is now completed. $\Box$

\bigskip
{\bf Proof of Corollary \ref{cor2}.} For $\lambda>0$, $\varphi_{\lambda}=(-\lambda)\vee (\varphi(x)\wedge \lambda)\in C_b(C[0,1])$. So, by Theorem \ref{thFCLT},
$$
\Sbep^{\ast}\left[\varphi_{\lambda}(W_n)\right]=\Sbep\left[\varphi_{\lambda}(W_n)\right]
\to \widetilde{\mathbb E}\left[\varphi_{\lambda}(W)\right]=\widetilde{\mathbb E}^{\ast}\left[\varphi_{\lambda}(W)\right].
$$
On the other hand,
\begin{align*}
& \left|\widetilde{\mathbb E}^{\ast}\left[\varphi_{\lambda}(W)\right]-\widetilde{\mathbb E}^{\ast}\left[\varphi(W)\right]\right|\le \widetilde{\mathbb E}^{\ast}\left[\left(|\varphi_{\lambda}(W)|-\lambda\right)^+\right]\\
& \qquad \le   C\widetilde{\mathbb E}^{\ast}\left[ \left(\|W\|^p-\lambda/C+1\right)^+\right]\to 0 \text{ as } \lambda\to 0
\end{align*}
and
\begin{align*}
 &    \left|\Sbep^{\ast}\left[\varphi_{\lambda}(W_n)\right]-\Sbep^{\ast}\left[\varphi(W_n)\right]\right| \le   C  \Sbep^{\ast}\left[ \left(\|W_n\|^p-\lambda/C+1\right)^+\right] \\
& \qquad   = C  \Sbep \left[ \left(\max_{k\le n}\left|\frac{S_k}{\sqrt{n}}\right|^p -\lambda/C+1\right)^+\right].
\end{align*}
It is sufficient to show
\begin{equation}\label{eqUnifIn}
\lim_{\lambda\to \infty} \limsup_n \Sbep\left[ \left(\max_{k\le n}\left|\frac{S_k}{\sqrt{n}}\right|^p -\lambda\right)^+\right]=0,
\end{equation}
i.e., the sequence $\{\max_{k\le n} \left|S_k/\sqrt{n}\right|^p;n\ge 1\}$ is uniformly integrable under $\Sbep$. Let   $Y_j=(-\sqrt{n})\vee(X_j\wedge \sqrt{n})$, $\widehat{X}_j=X_j-Y_j$, $T_j=\sum_{i=1}^j Y_i$, $\widehat{S}_j=\sum_{i=1}^j \widehat{X}_i$, $j=1,\ldots, n$.
Then $\max_{k\le n}|S_k|\le \max_{k\le n}|\widehat{S}_k|+\max_{k\le n}|T_k|$.
Note $\Sbep[X_1]=\cSbep[X_1]=0$. So, $|\cSbep[Y_1]|=|\cSbep[X_1]-\cSbep[Y_1]|\le \Sbep|\widehat{X}_1|=\Sbep[(|X_1|^2-n)^+]n^{-1/2}$ and
$|\Sbep[Y_1]|=|\Sbep[X_1]-\Sbep[Y_1]|\le \Sbep|\widehat{X}_1|=\Sbep[(|X_1|^2-n)^+]n^{-1/2}$.
By Rosnethal's inequality (c.f, (\ref{eqthRIneq.2})),
\begin{align*}
& \Sbep\left[\max_{k\le n} \left|T_k\right|^{2p}\right]\le    C_p\left\{ n \Sbep [|Y_1|^{2p}]+\left(n \Sbep [|Y_1|^2]\right)^{p}  +\left(n\big[\big(\cSbep [Y_1]\big)^-+\big(\Sbep [Y_1]\big)^+\big]\right)^{2p}\right\}\\
& \;\; \le  C_p\left\{ n n^{p/2}\Sbep [|X_1|^{p}]+n^p\left(\Sbep [X_1^2]\right)^{p}+\left(nn^{-1/2}\Sbep\left[(X_1^2-n)^+\right]\right)^{2p}\right\} \le C n^p
\end{align*}
and
\begin{align*}
 &\Sbep\left[\max_{k\le n} \left|\widehat{S}_k\right|^{p}\right]
\le    C_p\left\{ n \Sbep [|\widehat{X}_1|^{p}]+\left(n \Sbep [|\widehat{X}_1|^2]\right)^{\red{p/2}}  +\left(n\big[\big(\cSbep [\widehat{X}_1]\big)^-+\big(\Sbep [\widehat{X}_1]\big)^+\big]\right)^{p}\right\}\\
\le & C_p\left\{  n\Sbep \left[\big(|X_1|^p-n^{p/2}\big)^+\right]+\red{n^{p/2}\left(\Sbep \left[(X_1^2-n)^+\right]\right)^{p/2}}+n^{\red{p/2}}\left(\Sbep \left[(X_1^2-n)^+\right]\right)^{p}\right\}.
\end{align*}
It follows that
\begin{align*}
& \lim_{\lambda\to \infty} \limsup_n \Sbep\left[ \left(\max_{k\le n}\left|\frac{T_k}{\sqrt{n}}\right|^p -\lambda\right)^+\right] \\
\le & \lim_{\lambda\to \infty} \limsup_n \frac{1}{\lambda}\Sbep\left[\max_{k\le n}\left|\frac{T_k}{\sqrt{n}}\right|^{2p}\right]\le \lim_{\lambda\to \infty}  \frac{C}{\lambda}=0
\end{align*}
and
\begin{align*}
  \limsup_n  \Sbep\left[\max_{k\le n}\left|\frac{\widehat{S}_k}{\sqrt{n}}\right|^{p}\right]=0.
\end{align*}
The proof of (\ref{eqUnifIn}) is completed. $\Box$

\bigskip
{\bf Proof of Corollary \ref{cor3}.} We only give the proof of the first result. Let $\varphi(y)$ be a Lipschitz function such that $I\{y\ge x\}\le  \varphi(y)\le I\{y\ge x(1+\delta)\}$. Then by Theorem \ref{thFCLT},
\begin{align*}
 \limsup_{n\to \infty}& \Capc\left(\max_{k\le n}|S_k|/\sqrt{n}\ge x\right)=\limsup_{n\to \infty} \Capc\left(\sup_{0\le t\le 1}|W_n(t)|\ge x\right)\\
\le &  \limsup_{n\to \infty} \Sbep\left[\varphi\left(\sup_{0\le t\le 1}|W_n(t)|\right)\right]
= \widetilde{\mathbb E}\left[\varphi\left(\sup_{0\le t\le 1}|W(t)|\right)\right]\\
\le & \widetilde{\Capc}\left(\sup_{0\le t\le 1}|W(t)|\ge x(1+\delta)\right)
=P\left(\sup_{0\le t\le 1}\overline{\sigma}|B(t)|\ge x(1+\delta)\right),
\end{align*}
where the last inequality is due to (\ref{eqBrownian1}). Letting $\delta\to 0$ yields
$$\limsup_{n\to \infty} \Capc\left(\max_{k\le n}|S_k|/\sqrt{n}\ge x\right)
\le P\left(\sup_{0\le t\le 1}\overline{\sigma}|B(t)|\ge x\right). $$
By considering a function $\varphi(y)$ with $I\{y\ge x(1-\delta)\}\le  \varphi(y)\le I\{y\ge x\}$ instead, we can show that
$$\liminf_{n\to \infty} \Capc\left(\max_{k\le n}|S_k|/\sqrt{n}\ge x\right)
\ge P\left(\sup_{0\le t\le 1}\overline{\sigma}|B(t)|\ge x\right). $$
The proof is completed. $\Box$.

\bigskip
{\bf Proof of Theorem \ref{thChungLIL}.} The proof is based on small deviations for $\max_{i\le n}|S_i|$ which are proved in Section \ref{sectDeviation}. Let $0<\epsilon <1/2$,
$\beta(n)=\sqrt{\frac{n \pi^2}{8\log\log n}}$ and $x_n=(1+\epsilon)^{-1}\underline{\sigma} \beta(n)/\sqrt{n}$. Then by (\ref{eqDeviation0.1}) (c.f, Theorem \ref{thDeviation}),
$$\log \Capc\Big(\max_{i\le n}|S_i|\le \beta(n) (1+\epsilon)^{-1}\underline{\sigma}\Big)
\sim -\frac{\pi^2\underline{\sigma}^2}{8x_n^2}\sim -  (1+\epsilon)^2\log\log n. $$
Let $n_k=[e^{k/\log k}]$. Then $n_{k-1}/n_k\to 1$, $\beta(n_{k-1})/\beta(n_k)\to 1$ and $\log\log n_k\sim \log k$. It follows that
$$\sum_{k=1}^{\infty} \Capc\Big(\max_{i\le n_k}|S_i|\le \beta(n_k) (1+\epsilon)^{-1}\underline{\sigma}\Big)
\le C\sum_{k=1}^{\infty} k^{-(1+\epsilon) }<\infty. $$
Hence by the countable sub-additivity of $\Capc$,
\begin{align*}
&\Capc\Big(\max_{i\le n_k}|S_i|\le \beta(n_k) (1+\epsilon)^{-1}\underline{\sigma} \;\; i.o. \Big)\\
\le &\sum_{k=K}^{\infty} \Capc\Big(\max_{i\le n_k}|S_i|\le \beta(n_k) (1+\epsilon)^{-1}\underline{\sigma}\Big)
\to 0 \;\text{ as } K\to \infty.
\end{align*}
So,
$$\Capc\Big(\liminf_{k\to \infty}\frac{\max_{i\le n_k}|S_i|}{\beta(n_k)}\le (1+\epsilon)^{-1}\underline{\sigma}\Big)=0. $$
Note for $n_{k}\le n\le n_{k+1}$,
$$ \frac{\max_{i\le n}|S_i|}{\beta(n)}\ge \frac{\max_{i\le n_k}|S_i|}{\beta(n_k)}\cdot \frac{\beta(n_k)}{\beta(n_{k+1})}. $$
It follows that
$$\Capc\Big(\liminf_{n\to \infty}\frac{\max_{i\le n}|S_i|}{\beta(n)}\le (1+\epsilon)^{-1}\underline{\sigma}\Big)=0. $$
Note the continuity  of $\Capc$. Letting $\epsilon\to 0$ yields
\begin{equation}\label{eqproofCLIL.1}\Capc\Big(\liminf_{n\to \infty}\frac{\max_{i\le n}|S_i|}{\beta(n)}<\underline{\sigma}\Big)=0. \end{equation}

Next, we consider the lower bound. Let $n_k=[e^{k(\log k)^2}]$, then $n_{k-1}/n_k\sim e^{-(\log k)^2}\to 0$ and $\log\log n_k \sim \log k$. Let $x_{n_k}=(1-\epsilon)^{-1}\underline{\sigma} \beta(n_k)/\sqrt{n_k-n_{k-1}}$.
Then by (\ref{eqDeviation0.1}) (c.f, Theorem \ref{thDeviation}),
$$\log \Capc\Big(\max_{n_{k-1}<i\le n_k}|S_i\red{-S_{n_{k-1}}}|\le \beta(n_k) (1-\epsilon)^{-1}\underline{\sigma}\Big)
\sim -\frac{\pi^2\underline{\sigma}^2}{8x_{n_k}^2}\sim -  (1-\epsilon)^{-2}\log k. $$
So,
$$\sum_{k=1}^{\infty} \Capc\Big(\max_{n_{k-1}<i\le n_k}|S_i\red{-S_{n_{k-1}}}|\le \beta(n_k) (1-\epsilon)^{-1}\underline{\sigma}\Big)=\infty. $$
Let $\xi_k =\max_{n_{k-1}<i\le n_k}|S_i\red{-S_{n_{k-1}}}|- \beta(n_k) (1-\epsilon)^{-1}\underline{\sigma}$ and $\varphi(y)$ be a Lipschitz function such that $I\{y\le 0\}\le  \varphi(y)\le I\{y\le \epsilon\}$.
Then $\sum_{k=1}^{\infty}\Sbep[\varphi(\xi_k)]=\infty$. Note the independence and the continuity of $\cCapc$. We have
\begin{align*}  &\cCapc\left(\bigcap_{k=n}^{\infty}\{\xi_k>\epsilon\}\right)
= \lim_{N\to \infty}\cCapc\left(\bigcap_{k=n}^N\{\xi_k>\epsilon\}\right)
\le  \lim_{N\to \infty} \cSbep\left[\prod_{k=n}^N(1-\varphi(\xi_k))\right] \\
= & \lim_{N\to \infty} \prod_{k=n}^N\cSbep\left[(1-\varphi(\xi_k))\right]
\le  \lim_{N\to \infty} \prod_{k=n}^N\left(1-\Sbep[\varphi(\xi_k)]\right)
\le \exp\left\{-\sum_{k=n}^{\infty}\Sbep[\varphi(\xi_k)]\right\}=0.
\end{align*}
Hence
$$ \Capc\left( \{\xi_k\le \epsilon\} \; i.o.\right)=1-\cCapc\left(\bigcup_{n=1}^{\infty}\bigcap_{k=n}^{\infty}\{\xi_k>\epsilon\}\right)
=1-\lim_{n\to \infty}\cCapc\left(\bigcap_{k=n}^{\infty}\{\xi_k>\epsilon\}\right)=1, $$
i.e.,
 $$ \Capc\Big(\max_{n_{k-1}<i\le n_k}|S_i\red{-S_{n_{k-1}}}|\le \beta(n_k) (1-\epsilon)^{-1}\underline{\sigma}+\epsilon\;\; i.o. \Big)=1. $$
 It follows that
\begin{equation}\label{eqproofCLIL.2} \Capc\Big(\liminf_{k\to \infty}\frac{\max_{n_{k-1}<i\le n_k}|S_i\red{-S_{n_{k-1}}}|}{\beta(n_k)}\le (1-\epsilon)^{-1}\underline{\sigma}+\epsilon\Big)=1.
\end{equation}
On the other hand,
\begin{align*}
 \sum_{k=1}^{\infty}\Capc\Big( &  \max_{i\le n_{k-1}}|S_i| \ge \beta(n_k)/\sqrt{\log\log n_k}\Big)
\le   \sum_{k=1}^{\infty}\frac{\log\log n_k}{\beta^2(n_k)}\Sbep\left[\max_{i\le n_{k-1}}|S_i|^2\right]\\
\le &\sum_{k=1}^{\infty} C \frac{(\log\log n_k)^2}{n_k}n_{k-1} \Sbep[X_1^2]\le \sum_{k=1}^{\infty}C \frac{(\log k)^3}{ e^{(\log k)^2}}<\infty,
\end{align*}
where the second inequality is due to the  Rosenthal type inequality (c.f. Lemma \ref{thRIneq}). It follows that
 $$ \Capc\Big(  \limsup_{k \to \infty}\frac{\max_{i\le n_{k-1}}|S_i|}{\beta(n_k)}>0  \Big)=0, $$
 which, together with (\ref{eqproofCLIL.2}), implies
 \begin{equation}\label{eqproofCLIL.3} \Capc\Big(\liminf_{k\to \infty}\frac{\max_{i\le n_k}|S_i|}{\beta(n_k)}\le (1-\epsilon)^{-1}\underline{\sigma}+\epsilon \Big)=1.
\end{equation}
By the continuity of $\Capc$, letting $\epsilon\to 0$ yields
$$\Capc\Big(\liminf_{k\to \infty}\frac{\max_{i\le n_k}|S_i|}{\beta(n_k)}\le  \underline{\sigma} \Big)=1.$$
Hence
\begin{equation}\label{eqproofCLIL.4} \Capc\Big(\liminf_{n\to \infty}\frac{\max_{i\le n}|S_i|}{\beta(n)}\le  \underline{\sigma} \Big)=1.
\end{equation}

Finally, let $x_{n_k}=(1-\epsilon)^{-1}\overline{\sigma} \beta(n_k)/\sqrt{n_k-n_{k-1}}$.
Then by (\ref{eqDeviation0.2}) (c.f, Theorem \ref{thDeviation}),
$$\log \cCapc\Big(\max_{n_{k-1}<i\le n_k}|S_i\red{-S_{n_{k-1}}}|\le \beta(n_k) (1-\epsilon)^{-1}\overline{\sigma}\Big)
\sim -\frac{\pi^2\overline{\sigma}^2}{8x_{n_k}^2}\sim -  (1-\epsilon)^2\log k, $$
which implies
$$\sum_{k=1}^{\infty} \cCapc\Big(\max_{n_{k-1}<i\le n_k}|S_i\red{-S_{n_{k-1}}}|\le \beta(n_k) (1-\epsilon)^{-1}\overline{\sigma}\Big)=\infty. $$
So, similar to (\ref{eqproofCLIL.4}) we have
\begin{equation}\label{eqproofCLIL.5} \cCapc\Big(\liminf_{n\to \infty}\frac{\max_{i\le n}|S_i|}{\beta(n)}\le  \overline{\sigma} \Big)=1.
\end{equation}
It is obvious that (\ref{eqChungLIL.1}) follows from  (\ref{eqproofCLIL.1}), (\ref{eqproofCLIL.5})
and the fact  $\cCapc(A\bigcap B)\ge \cCapc(A)-\Capc(B^c)$, and
(\ref{eqChungLIL.2}) follows from  (\ref{eqproofCLIL.1}), (\ref{eqproofCLIL.4}) and the fact  $\Capc(A\bigcap B)\ge \Capc(A)-\Capc(B^c)$. $\Box$

\section{Convergence of the finite-dimensional distributions}\label{sectFinite}
\setcounter{equation}{0}

The purpose of this section is to prove the following theorem on the convergence of finite-dimensional distributions of $W_n$ under $\Sbep$.
\begin{theorem} \label{thfinite}    Suppose $\Sbep[(X_1^2-b)^+]\to 0$ as $b\to \infty$. Then for any $0\le t_1<\ldots<t_d\le 1$ and any $\varphi\in C_b(\mathbb R^d)$,   we have
\begin{equation}
\Sbep\left[\varphi\left(W_n(t_1),\ldots, W_n(t_d)\right)\right]\to \widetilde{\mathbb E}\left[\varphi\left(W(t_1),\ldots, W(t_d)\right)\right].
\end{equation}
\end{theorem}

\bigskip
For proving this theorem, we need some lemmas. The first is the central limit theorem which was firstly obtained  by Peng (2008b) under the condition $\Sbep[|X_1|^{2+\epsilon}]<\infty$.  The condition  was relaxed to only the existence of the second moments by Zhang (2014).

 \begin{lemma}\label{thCLT} (CLT) Suppose $\Sbep[(X_1^2-b)^+]\to 0$ as $b\to \infty$. Then for any  bounded continuous function $\varphi$,
  \begin{equation} \label{eqCLT} \lim_{n\to \infty}\Sbep\left[\varphi\left(\frac{S_n}{\sqrt{n}}\right)\right]=\widetilde{\mathbb E}[\varphi(W(1))].
  \end{equation}
 \end{lemma}

For random vectors $\bm X_n=(X_n^1,\cdots, X_n^d)$ in $\mathscr{H}^d$ and $\bm X=(X^1,\cdots, X^d)$ in $\widetilde{\mathscr{H}}^d$, we write $\bm X_n \overset{d}\to \bm X$ if
$ \Sbep[\varphi(\bm X_n)]\to \widetilde{\mathbb E}[\varphi(\bm X)]$
for any bounded continuous function
$\varphi: \mathbb R^d\to \mathbb R$;
and write $\bm X_n\overset{\Capc }\to \bm a$ if
$  \Capc\left(\|\bm X_n-\bm a\|>\epsilon\right)\to 0 \text{ for all } \epsilon>0. $ It is obvious that for any continuous function $f(\bm x)$, $\bm X_n \overset{d}\to \bm X$  implies $f(\bm X_n) \overset{d}\to f(\bm X)$, and $\bm X_n \overset{\Capc}\to \bm X$  implies $f(\bm X_n) \overset{\Capc}\to f(\bm X)$.

The following lemma is  Slutsky's theorem. The proof is standard and ommited.
\begin{lemma}\label{lemSlutsky} Suppose $\bm X_n\overset{d}\to \bm X$, $\bm Y_n \overset{\Capc}\to \bm y$, $\eta_n\overset{\Capc}\to a$, where $a$ is a constant and $\bm y$ is a constant vector, and $\widetilde{\Capc}(\|\bm X\|>\lambda)\to 0$ as $\lambda\to \infty$. Then
$(\bm X_n, \bm Y_n, \eta_n)\overset{d}\to (\bm X,\bm y, a)$, and as a result,
$\eta_n\bm X_n+\bm Y_n\overset{d}\to a\bm X+\bm y$.
\end{lemma}

\begin{lemma}\label{lemfinite3} Suppose $\bm X_n\overset{d}\to \bm X$ and $\widetilde{\Capc}(\|\bm X\|>\lambda)\to 0$ as $\lambda\to \infty$. Assume that $g_n(\bm x)$ and $g(x)$ are  continuous functions for which
$$ |g_n(\bm x)|\le M, \;\; |g(\bm x)|\le M, $$
$$ g_n(\bm x_n)\to g(\bm x) \text{ whenever } x_n \to x. $$
Then
$$\Sbep[g_n(\bm X_n)]\to \widetilde{\mathbb E}[g(\bm X)]. $$
\end{lemma}
{\bf Proof.} The conditions for $g_n$ imply that
$ \sup_{\|\bm x\|\le \lambda}|g_n(\bm x)-g(\bm x)|\to 0. $
It is obvious that
\begin{align*}
&\left|\Sbep[g_n(\bm X_n)]-\widetilde{\mathbb E}[g(\bm X)]\right|
\le \left|\Sbep[g_n(\bm X_n)]-\Sbep[g(\bm X_n)]\right|
+\left|\Sbep[g(\bm X_n)]-\widetilde{\mathbb E}[g(\bm X)]\right|\\
\le & \left|\Sbep[g(\bm X_n)]-\widetilde{\mathbb E}[g(\bm X)]\right|+\sup_{\|\bm x\|\le \lambda} |g_n(\bm x)-g(\bm x)|
+2 M \Capc\left(\|\bm X_n\|>\lambda\right).
\end{align*}
Choose a Lipschitz   function $\varphi(x)$ such that $I\{x>\lambda\}\ge \varphi(x)\ge I\{x>\lambda/2\}$. Letting
$n\to \infty$ yields that
\begin{align*}
&\limsup_{n\to \infty}\left|\Sbep[g_n(\bm X_n)]-\widetilde{\mathbb E}[g(\bm X)]\right|\\
\le & 2 M \limsup_{n\to \infty} \Capc\left(\|\bm X_n\|>\lambda\right)
\le 2 M \limsup_{n\to \infty}  \Sbep\left[\varphi(\|\bm X_n\|)\right]\\
=& 2 M \widetilde{\mathbb E}\left[\varphi(\|\bm X\|)\right]\le 2 M \widetilde{\Capc}\left(\|\bm X\|>\lambda/2\right)\to \text{ as } \lambda\to \infty.
\end{align*}
The proof is completed. $\Box$

\begin{lemma}\label{lemfinite4} Suppose that $\bm X_n\overset{d}\to \bm X$, $\bm Y_n\overset{d}\to \bm Y$,
 $\bm Y_n$ is independent to $\bm X_n$ under $\Sbep$, $\bm Y$ is independent to $\bm X$ under $\widetilde{\mathbb E}$, and $\widetilde{\Capc}(\|\bm X\|>\lambda)\to 0$  and $\widetilde{\Capc}(\|\bm Y\|>\lambda)\to 0$  as $\lambda\to \infty$. Then
 $ (\bm X_n,\bm Y_n)\overset{d }\to  (\bm X,\bm Y). $
\end{lemma}
{\bf Proof.} Suppose $\varphi(\bm x,\bm y)$ is a bounded continuous function. We want to show that
\begin{equation}\label{eqlemfinite4.1} \Sbep\left[\varphi(\bm X_n,\bm Y_n)\right]\to \widetilde{\mathbb E}\left[\varphi(\bm X,\bm Y)\right].
\end{equation}
 First we assume that $\varphi(\bm x,\bm y)$ is a bounded Lipschitz function. Then $\varphi\in C_{l,Lip}$. By the definition of the independence,
 $$ \Sbep\left[\varphi(\bm X_n,\bm Y_n)\right]=\Sbep\left[g_n(\bm X_n)\right], \;\;
 \widetilde{\mathbb E}\left[\varphi(\bm X,\bm Y)\right]=\widetilde{\mathbb E}\left[g(\bm X)\right], $$
where
$ g_n(\bm x)=\Sbep[\varphi(\bm x, \bm Y_n)]$, $ g(\bm x)=\widetilde{\mathbb E}[\varphi(\bm x, \bm Y)]$.
Suppose $\bm x_n\to \bm x$.   It follows that
\begin{align*}
& \left| g_n(\bm x_n)-g(\bm x)\right|=\left|\Sbep[\varphi(\bm x_n, \bm Y_n)]- \widetilde{\mathbb E} [\varphi(\bm x, \bm Y)]\right| \\
 \le & \left|\Sbep[\varphi(\bm x_n, \bm Y_n)]- \Sbep[\varphi(\bm x, \bm Y_n)]\right|
 +\left|\Sbep[\varphi(\bm x, \bm Y_n)]- \widetilde{\mathbb E}[\varphi(\bm x, \bm Y)]\right| \\
 \le &  \sup_{\bm y} \left| \varphi(\bm x_n, \bm y) -  \varphi(\bm x, \bm y)\right|
  +\left|\Sbep[\varphi(\bm x, \bm Y_n)]- \widetilde{\mathbb E}[\varphi(\bm x, \bm Y)]\right|\to 0
 \end{align*}
 by noting that $\varphi(\bm x,\bm y)$ is uniformly continuous  and $\bm Y_n\overset{d}\to \bm Y$.
 By the uniform continuity of $\varphi$, $g_n(\bm x)$ and $g(\bm x)$ are continuous functions.
So, $g_n(\bm x)$ and $g(\bm x)$ satisfy the conditions in Lemma \ref{lemfinite3}. It follows that
\begin{align*}
\Sbep[\varphi(\bm X_n,\bm Y_n)]=\Sbep[g_n(\bm X_n)]\to \widetilde{\mathbb E}[g(\bm X)]=
 \widetilde{\mathbb E}[\varphi(\bm X,\bm Y)].
\end{align*}
Next, we assume that $\varphi(\bm x,\bm y)$ is a bounded uniformly continuous function. Then for any $\epsilon>0$, there is bounded Lipschitz function $\varphi_{\epsilon}(\bm x,\bm y)$ such that $|\varphi(\bm x,\bm y)-\varphi_{\epsilon}(\bm x,\bm y)|<\epsilon$. It follows that
\begin{align*}
& \limsup_{n\to \infty} \left|\Sbep[\varphi(\bm X_n,\bm Y_n)]-\widetilde{\mathbb E}[\varphi(\bm X,\bm Y)\right| \\
\le & \limsup_{n\to \infty} \left|\Sbep[\varphi_{\epsilon}(\bm X_n,\bm Y_n)]-\widetilde{\mathbb E}[\varphi_{\epsilon}(\bm X,\bm Y)]\right|+2\epsilon=2\epsilon.
\end{align*}
So, (\ref{eqlemfinite4.1}) is proved for a bounded uniformly continuous function. Finally, let $\varphi(\bm x,\bm y)$ be a  bounded continuous function with $|\varphi(\bm x,\bm y)|\le M$. Let $\lambda>0$. For $\bm x=(x_1,\ldots, x_d)$, denote $\bm x_{\lambda}=\big((-\lambda)\vee(x_1\wedge \lambda)\lambda,\ldots, (-\lambda)\vee(x_d\wedge \lambda)\big)$ and define $\varphi_{\lambda}(\bm x,\bm y)=\varphi(\bm x_{\lambda},\bm y_{\lambda})$. Then $\varphi_{\lambda}$ is a bounded uniformly continuous function with
$$|\varphi_{\lambda}(\bm x,\bm y)-\varphi(\bm x,\bm y)|\le 2MI\{\|\bm x\|>\lambda\}+2MI\{\|\bm y\|>\lambda\}. $$
It follows that
\begin{align*}
& \limsup_{n\to \infty} \left|\Sbep[\varphi(\bm X_n,\bm Y_n)]-\widetilde{\mathbb E}[\varphi(\bm X,\bm Y)\right| \\
\le & \limsup_{n\to \infty} \left|\Sbep[\varphi_{\lambda}(\bm X_n,\bm Y_n)]-\widetilde{\mathbb E}[\varphi_{\lambda}(\bm X,\bm Y)]\right| \\
& + 2M\limsup_{n\to \infty}
\big\{\Capc(\|\bm X_n\|>\lambda)+\Capc(\|\bm Y_n\|>\lambda)\big\}\\
& +
2M
\big\{\widetilde{\Capc}(\|\bm X\|>\lambda)+\widetilde{\Capc}(\|\bm Y\|>\lambda)\big\}\\
\le & 4M\big\{\widetilde{\Capc}(\|\bm X\|>\lambda/2)+\widetilde{\Capc}(\|\bm Y\|>\lambda/2)\big\}\to 0\; \text{ as } \lambda\to \infty.
\end{align*}
 The proof is completed. $\Box$

\bigskip
{\bf Proof of Theorem \ref{thfinite}}. Note
$$ \sup_{0\le t\le 1}\left|W_n(t)-\frac{S_{[nt]}}{\sqrt{n}}\right|\le \frac{\max_{k\le n}|X_k|}{\sqrt{n}}, $$
and for any $\epsilon>0$,
$$ \Capc\left(\max_{k\le n}|X_k|\ge \epsilon \sqrt{n}\right)
\le n \frac{2}{\epsilon^2n}\Sbep\left[\big(X_1^2-\frac{\epsilon^2 n}{2}\big)^+\right]\to 0. $$
It is sufficient to show that
$$ \frac{1}{\sqrt{n}}\left( S_{[nt_1]},\ldots, S_{[nt_d]}\right)\overset{d}
\to \left( W(t_1),\ldots, W(t_d)\right)$$
by Lemma \ref{lemSlutsky},  or equivalently,
\begin{align}\label{eqproofofthfinite}
 &\frac{1}{\sqrt{n}}\left( S_{[nt_1]}, S_{[nt_2]}-S_{[nt_1]},\ldots, S_{[nt_d]}-S_{[nt_{d-1}]}\right)\nonumber  \\
& \quad \overset{d}\to \left( W(t_1), W(t_2)-W(t_1),\ldots, W(t_d)-W(t_{d-1})\right).
\end{align}
By Lemmas \ref{thCLT} and \ref{lemfinite3},
$$ \frac{S_{[nt_i]}-S_{[nt_{i-1}]}}{\sqrt{n}}=\frac{\sqrt{[nt_i]-[nt_{i-1}]}}{\sqrt{n}}
\frac{S_{[nt_i]}-S_{[nt_{i-1}]}}{\sqrt{[nt_i]-[nt_{i-1}]}}\overset{d}\to W(t_i)-W(t_{i-1}). $$
Hence, by noting the independence, (\ref{eqproofofthfinite}) follows from Lemma \ref{lemfinite4} and the induction.
The proof is now completed. $\Box$

\section{Tightness}\label{sectTight}
\setcounter{equation}{0}

Recall $ \omega_{\delta}(x)=\sup\limits_{|t-s|<\delta}|x(t)-x(s)|$. It is known that $\{W_n\}$ is tight under the probability measure $P$ in the following sense: for any $\eta>0$, there a compact set $K\subset C[0,1]$ such that $\sup_n P(W_n\not\in K)<\eta$. This is also equivalent to $\lim_{\delta\to 0} \sup_nP(\omega_{\delta}(W_n)\ge \epsilon)=0$ for any $\epsilon>0$.
In this section, we will prove the following theorem  on the tightness of $\{W_n\}$ under capacities.
\begin{theorem}\label{thTight} Suppose $\Sbep[(X_1^2-b)^+]\to 0$ as $b\to \infty$. Then
\begin{description}
  \item[\rm (a)] for any $\epsilon>0$,
\begin{equation}\label{eqTight.1}\lim_{\delta\to 0}\sup_n \Capc\left(\omega_{\delta}(W_n)\ge \epsilon\right)=0; \end{equation}
  \item[\rm (b)] for any $\eta>0$,  there exists a compact set $K\subset C[0,1]$ such that
 \begin{equation}\label{eqTight.2} \sup_n \Capc\left(W_n\not\in K\right)<\eta.
 \end{equation}
\end{description}
\end{theorem}

 {\bf Proof.} We first show (a).
 With the same argument of Billingsley (1968, Pages 56-59, c.f., Theorem 8.4), it is sufficient to show that
 \begin{equation}\label{eqProofTight.2}
\lim_{\lambda\to \infty}\lambda^2\sup_n \sup_k\Capc\left(\max_{i\le n}|S_{k+i}-S_k|\ge \lambda \sqrt{n}\right)=0.
 \end{equation}
Note that for each fixed $n$,
\begin{align*}
&\lim_{\lambda\to \infty}\lambda^2 \sup_k\Capc\left(\max_{i\le n}|S_{k+i}-S_k|\ge \lambda \sqrt{n}\right)\\
\le & \sum_{i=1}^n \lim_{\lambda\to \infty}\sup_k\Capc\left(|X_{k+i}|\ge \lambda/\sqrt{n}\right)
\le  2n^2 \lim_{\lambda\to \infty}\Sbep\left[ \big(X_1^2-\frac{\lambda^2}{2n}\big)^+\right]=0.
\end{align*}
So, it is sufficient to show that
 \begin{equation}\label{eqProofTight.2}
\lim_{\lambda\to \infty}\lambda^2\limsup_{n\to \infty} \sup_k\Capc\left(\max_{i\le n}|S_{k+i}-S_k|\ge \lambda \sqrt{n}\right)=0.
 \end{equation}
 Now,
 \begin{align*}
  \Capc\left(\max_{i\le n}|S_{k+i}-S_k|\ge \lambda \sqrt{n}\right)=  \Capc\left(\max_{i\le n}|S_i|\ge \lambda \sqrt{n}\right)
 \le   \frac{2}{\lambda^2}\Sbep\left[\left(\max_{i\le n}\frac{|S_i|^2}{n}-\frac{\lambda^2}{2}\right)^+\right].
 \end{align*}
By (\ref{eqUnifIn}) where $p=2$, (\ref{eqProofTight.2}) follows.

 Now, we show (\ref{eqTight.2}). By (\ref{eqTight.1}), choose $\delta_k\downarrow 0$ such that, if
 $$A_k=\left\{x: \omega_{\delta_k}(x)<\frac{1}{k}\right\}, $$
then $\sup_n\Capc\left( W_n \in A_k^c\right)\le \eta/2^{k+1}$.
Let $A=\{x:|x(0)|\le a\}$, $K=A\bigcap_{k=1}^{\infty}A_k$. Then by the Arzel\'a-Ascoli theorem, $K$ is compact.  It is obvious that
$\{W_n\not\in A\}=\emptyset$ since $W_n(0)=0$. Next, we show that
$$\Capc(W_n\in K^c)\le \Capc(W_n\in A^c)+\sum_{k=1}^{\infty}\Capc(W_n\in A_k^c), $$
which is obvious when $\Capc$ is countably  sub-additive. Note that when $\delta<1/(2n)$,
$$ \omega_{\delta}(W_n)\le 2 |t-s|\max_{i\le n} \frac{|X_i|/\sqrt{n}}{1/n}\le 2 \sqrt{n} \delta  \max_{i\le n} |X_i|. $$
Choose a $k_0$ such that $\delta_k<1/(2Mk)$ and $\delta_k<1/(2n)$ for $k\ge k_0$. Then on the event $ E=\{\sqrt{n}\max_{i\le n} |X_i|\le M\}$, $\{W_n\in A_k^c\}=\emptyset$ for $k\ge k_0$. So,
by the (finite) sub-additivity of $\Capc$,
\begin{align*}
\Capc\big(E \bigcap \{W_n\in K^c\}\big)\le & \Capc\big(E \bigcap\{W_n \in  A^c\}\big)+\sum_{k=1}^{k_0}\Capc\big(E\bigcap \{W_n \in A_k^c\}\big) \\
\le & \Capc(W_n \in A^c)+\sum_{k=1}^{\infty}\Capc(W_n \in A_k^c).
\end{align*}
On the other hand,
$$\Capc(E^c )\le \frac{\sqrt{n}\Sbep[\max_{i\le n} |X_i|]}{M}\le \frac{n^2 \Sbep[|X_1|]}{M}. $$
It follows that
\begin{align*}
\Capc\big( W_n\in K^c \big)
\le  \Capc(W_n \in A^c)+\sum_{k=1}^{\infty}\Capc(W_n \in A_k^c)+ \frac{n^2 \Sbep[|X_1|]}{M}.
\end{align*}
Letting $M\to \infty$ yields
\begin{align*} \Capc\big( W_n\in K^c \big)\le
   \Capc(W_n \in A^c)+\sum_{k=1}^{\infty}\Capc(W_n \in A_k^c)
<   0+\sum_{k=1}^{\infty} \frac{\eta}{2^{k+1}}<\eta.
\end{align*}
The proof of (\ref{eqTight.2}) is now completed. $\Box$.

For the G-Brownian  motion $W(t)$ on $(\widetilde{\Omega},\widetilde{\mathscr{H}},\widetilde{\mathbb E})$ we have a similar result.
 \begin{theorem}\label{thTight2}   We have
\begin{description}
  \item[\rm (a)] for any $\epsilon>0$,
$\lim_{\delta\to 0} \widetilde{\Capc}\left(\omega_{\delta}(W)\ge \epsilon\right)=0; $
  \item[\rm (b)] for any $\eta>0$,  there exists a compact set $K\subset C[0,1]$ such that
$ \widetilde{\Capc}\left(W\not\in K\right)<\eta.$
\end{description}
\end{theorem}

{\bf Proof.} Note
\begin{align*}
 \widetilde{\Capc}\left(\omega_{\delta}(W)\ge \epsilon\right)
\le & \sum_{k=0}^{[1/\delta]} \widetilde{\Capc}\left(\sup_{0\le s\le \delta} |W(k/\delta+s)-W(k/\delta)|\ge \epsilon\right)\\
\le & \sum_{k=0}^{[1/\delta]}\frac{1}{\epsilon^4}\widetilde{\mathbb E}\left[\sup_{0\le s\le \delta} |W(k/\delta+s)-W(k/\delta)|^4 \right] \\
\le & \frac{2}{\delta\epsilon^4} \delta^2\widetilde{\mathbb E}\left[\sup_{0\le s\le 1} | W(s)|^4 \right]
=  \delta \frac{2\overline{\sigma}^4}{\epsilon^4} E_P \left[\sup_{0\le s\le 1} | B(s)|^4 \right].
\end{align*}
Hence (a) follows. The proof of (b) is similar to that of Theorem \ref{thTight} (b) by noting that $ \widetilde{\Capc}$ is countably sub-additive. $\Box$

\bigskip
In the end of this section, we give another proof of Theorem \ref{thFCLT}. Define
$$ {\mathbb F}_n[\varphi]=\Sbep\left[\varphi(W_n)\right], \;  {\mathbb F}[\varphi]=\sup_n\Sbep\left[\varphi(W_n)\right], \; \; \varphi\in C_b\left(C[0,1]\right) $$
and ${\mathbb F}_n^{\ast}[\varphi]$, ${\mathbb F}^{\ast}$ be their extensions. By Theorem \ref{thTight}, $\mathbb F$ is tight and hence the family of sub-linear expectations $\{{\mathbb F}_n; n\ge 1\}$ on $\left(C[0,1],C_b\left(C[0,1]\right)\right)$  is tight in the sense of Definition 7 of Peng (2010). Hence, by Theorem 9 of Peng (2010), $\{{\mathbb F}_n\}$ is weakly compact, namely, for each subsequence $\{{\mathbb F}_{n^{\prime}}\}$ there
exists a further subsequence $\{{\mathbb F}_{n^{\prime\prime}}\}$ such that, for each $\varphi\in C_b\left(C[0,1]\right)$, $\{{\mathbb F}_{n^{\prime\prime}}[\varphi]\}$  is
a Cauchy sequence. Define
$$ \widetilde{\mathbb F}[\varphi]=\lim_{n^{\prime\prime}\to \infty}{\mathbb F}_{n^{\prime\prime}}[\varphi]
,\; \varphi\in C_b\left(C[0,1]\right). $$
Then $\widetilde{\mathbb F}$ satisfies (\ref{eqBM}) by Theorem \ref{thfinite}. So, under the sub-linear expectation $\widetilde{\mathbb F}$, the canonical process $W(t)=\omega_t$ is a G-Browinian with $W(1)\sim N(0,[\underline{\sigma}^2, \overline{\sigma}^2])$. The proof is completed.

\section{Small deviations}\label{sectDeviation}
\setcounter{equation}{0}

The purpose of this section is to establish the following theorem on the small deviations under $\Sbep$.
\begin{theorem}\label{thDeviation} Suppose $\Sbep[(X_1^2-b)^+]\to 0$ as $b\to \infty$, $0<x_n\to 0$ and $n^{1/2}x_n\to \infty$. Then
\begin{align}\label{eqDeviation.1}
&\lim_{n\to \infty} x_n^2\log \Capc\left(\max_{k\le n}|S_k|\le n^{1/2}x_n\right)=-\frac{\pi^2 \underline{\sigma}^2}{8},\\
\label{eqDeviation.2}
&\lim_{n\to \infty} x_n^2\log \cCapc\left(\max_{k\le n}|S_k|\le n^{1/2}x_n\right)=-\frac{\pi^2 \overline{\sigma}^2}{8}.
\end{align}
\end{theorem}

To prove Theorem \ref{thDeviation}, we need some lemmas on the properties of G-Brownian motions.
\begin{lemma}\label{lemma7} We have for all $x>0$,
\begin{align}\label{eqlemma7.1}
&\widetilde{\Capc}\left(\sup_{0\le t\le 1}|W(t)|\le x\right)=\sup_{\underline{\sigma}\le \sigma
\le \overline{\sigma}}P\left(\sup_{0\le t\le 1}|\sigma B(t)|\le x\right)
=P\left(\sup_{0\le t\le 1}|\underline{\sigma} B(t)|\le x\right)\\
\label{eqlemma7.2}
&\widetilde{\cCapc}\left(\sup_{0\le t\le 1}|W(t)|\le x\right)=\inf_{\underline{\sigma}\le \sigma
\le \overline{\sigma}}P\left(\sup_{0\le t\le 1}|\sigma B(t)|\le x\right)
=P\left(\sup_{0\le t\le 1}|\overline{\sigma} B(t)|\le x\right) \\
\label{eqlemma7.3}
&\widetilde{\Capc}\left(\sup_{0\le t\le 1} W(t) \le x\right)=\sup_{\underline{\sigma}\le \sigma
\le \overline{\sigma}}P\left(\sup_{0\le t\le 1} \sigma B(t) \le x\right)
=P\left(\sup_{0\le t\le 1} \underline{\sigma} B(t) \le x\right)\\
\label{eqlemma7.4}
&\widetilde{\cCapc}\left(\sup_{0\le t\le 1} W(t) \le x\right)=\inf_{\underline{\sigma}\le \sigma
\le \overline{\sigma}}P\left(\sup_{0\le t\le 1} \sigma B(t) \ge x\right)
=P\left(\sup_{0\le t\le 1} \overline{\sigma} B(t) \le x\right).
\end{align}
In particular, for $x>0$,
\begin{align}\label{eqlemma7.5}
&\frac{2}{\pi}\exp\Big\{-\frac{\pi^2 \underline{\sigma}^2}{8x^2}\Big\}
\le \widetilde{\Capc}\left(\sup_{0\le t\le 1}|W(t)|\le x\right)\le \frac{4}{\pi}\exp\Big\{-\frac{\pi^2 \underline{\sigma}^2}{8x^2}\Big\},\\
\label{eqlemma7.6}
&\frac{2}{\pi}\exp\Big\{-\frac{\pi^2 \overline{\sigma}^2}{8x^2}\Big\}
\le \widetilde{\cCapc}\left(\sup_{0\le t\le 1}|W(t)|\le x\right)\le \frac{4}{\pi}\exp\Big\{-\frac{\pi^2 \overline{\sigma}^2}{8x^2}\Big\}.
\end{align}
\end{lemma}
{\bf Proof.} Let $\varphi(y)$ be a non-increasing Lipschitz function for which $I\{y\le x\}\le \varphi(y)\le I\{y\le x(1+\delta)\}$. Then by Lemma \ref{DenisHuPeng},
\begin{align*}
&\widetilde{\Capc}\left(\sup_{0\le t\le 1}|W(t)|\le x\right)\le \widetilde{\mathbb E}\left[\varphi\big(\sup_{0\le t\le 1}|W(t)|\big)\right] \\
= & \sup_{\theta\in \Theta} E_P\left[\varphi\left(\sup_{0\le t\le 1}\Big|\int_0^t\theta(s)d B(s)\Big|\right)\right]\le
\sup_{\theta\in \Theta}P\left( \sup_{0\le t\le 1}\Big|\int_0^t\theta(s)d B(s)\Big|\le x(1+\delta)\right).
\end{align*}
Note that $W_{\theta}(t)=\int_0^t\theta(s)d B(s)$ is a continuous martingale with quadratic variation process $\langle W_{\theta},W_{\theta}\rangle (t)=\int_0^t\theta^2(s)ds$. By the Dambis-Dubins-Schwarz theorem, there is a standard Brownian motion $\widetilde{B}$ under $P$ such  that   $W_{\theta}(t)=\widetilde{B}\left(\langle W_{\theta},W_{\theta} \rangle(t)\right)$. On the other hand, $\langle W_{\theta},W_{\theta}\rangle(t)$ is continuous.  So,
$$\sup_{0\le t\le \overline{\sigma}^2 }|\widetilde{B}(s)|\ge \sup_{0\le t\le 1}\left|\widetilde{B}\left(\langle W_{\theta},W_{\theta}\rangle(t)\right)\right|
=\sup_{0\le t\le \langle W_{\theta},W_{\theta}\rangle(1)}|\widetilde{B}(s)|\ge \sup_{0\le t\le \underline{\sigma}^2 }|\widetilde{B}(s)|. $$
It follows that
$$\widetilde{\Capc}\Big(\sup_{0\le t\le 1}|W(t)|\le x\Big)\le
P\Big( \sup_{0\le t\le \underline{\sigma}^2 }|B(s)|\le x(1+\delta)\Big). $$
Letting $\delta\to 0$ yields
 $$\widetilde{\Capc}\Big(\sup_{0\le t\le 1}|W(t)|\le x\Big)\le
P\Big( \sup_{0\le t\le \underline{\sigma}^2 }|\widetilde{B}(s)|\le x\Big)=P\Big( \sup_{0\le t\le 1} \underline{\sigma}|B(s)|\le x\Big). $$
On the other hand, for $\theta(s)\equiv \underline{\sigma}$, $W_{\theta}(t)$ is a Brownian motion under $P$  with $W(1)\sim N(0,\underline{\sigma}^2)$. So
$$\widetilde{\Capc}\Big(\sup_{0\le t\le 1}|W(t)|\le x\Big)\ge
P\Big( \sup_{0\le t\le \underline{\sigma}^2 }|B(s)|\le x\Big). $$
It follows that
$$\widetilde{\Capc}\Big(\sup_{0\le t\le 1}|W(t)|\le x\Big) =P\Big( \sup_{0\le t\le 1} \underline{\sigma}|B(s)|\le x\Big). $$
Hence, (\ref{eqlemma7.1}) is proved. The proof of (\ref{eqlemma7.2})-(\ref{eqlemma7.4}) is similar.
The proof is completed by noting
$$ \frac{2}{\pi}\exp\Big\{-\frac{\pi^2  }{8x^2}\Big\}
\le P\Big(\sup_{0\le t\le 1}|B(t)|\le x\Big)\le \frac{4}{\pi}\exp\Big\{-\frac{\pi^2 }{8x^2}\Big\}.\;\; \Box $$

\begin{lemma}\label{lemma8} We have for all $y$,
\begin{align}\label{eqLemma8.1}
&  \widetilde{\Capc}\left(\sup_{0\le t\le 1}|W(t)+y |\le x\right)\le \frac{4}{\pi}\exp\Big\{-\frac{\pi^2 \underline{\sigma}^2}{8x^2}\Big\}\\
\label{eqLemma8.2}
&  \widetilde{\cCapc}\left(\sup_{0\le t\le 1 }|W(t)+y|\le x\right)\le \frac{4}{\pi}\exp\Big\{-\frac{\pi^2 \overline{\sigma}^2}{8x^2}\Big\}.
\end{align}
\end{lemma}
{\bf Proof.} The proof is similar to that of Lemma \ref{lemma7} by noting
$$P\Big(\sup_{0\le t\le 1}|B(t)+y|\le x\Big)\le P\Big(\sup_{0\le t\le 1}|B(t)|\le x\Big)\le \frac{4}{\pi}\exp\Big\{-\frac{\pi^2 }{8x^2}\Big\}  $$
according to the Anderson inequality.  $\Box$

\begin{lemma}\label{lemma9} We have for any $\alpha>0$, $0<\epsilon<\alpha/2$ and $\delta>0$,
\begin{align}\label{eqLemma9.1}
&\liminf_{x\to 0^+} x^2\log\left[\inf_{|y|\le \epsilon x}\Capc\Big(\sup_{0\le t\le 1}|W(t)|\le \alpha x,
     |y+W(1)|\le \delta x\Big)\right]\ge -\frac{\pi^2\underline{\sigma}^2}{8} (\alpha-2\epsilon)^{-2},\\
     \label{eqLemma9.2}
&\liminf_{x\to 0^+} x^2\log\left[\inf_{|y|\le \epsilon x}\cCapc\Big(\sup_{0\le t\le 1}|W(t)|\le \alpha x,
     |y+W(1)|\le \delta x\Big)\right]\ge -\frac{\pi^2\overline{\sigma}^2}{8} (\alpha-2\epsilon)^{-2},
\end{align}
\end{lemma}
{\bf Proof.} Note
$$\Capc\Big(\sup_{0\le t\le 1}|W(t)|\le \alpha x,
     |y+W(1)|\le \delta x\Big)\ge P \Big(\sup_{0\le t\le 1}|B(t\underline{\sigma}^2)|\le \alpha x,
     |y+W(\underline{\sigma}^2)|\le \delta x\Big). $$
By Lemma 3.2 of Acosta (1983),
$$\liminf_{x\to 0^+} x^2\log\left[\inf_{|y|\le \epsilon x}P \Big(\sup_{0\le t\le 1}|B(t\underline{\sigma}^2)|\le \alpha x,
     |y+W(\underline{\sigma}^2)|\le \delta x\Big)\right]\ge -\frac{\pi^2\underline{\sigma}^2}{8} (\alpha-2\epsilon)^{-2}. $$
The proof of (\ref{eqLemma9.1}) is completed.

The proof of (\ref{eqLemma9.2}) is more technical and similar to that of Lemma 3.2 of Acosta (1983) after smoothing  $I\{y\le r\}$ by a Lipschitz function $\varphi(y)$, and so omitted.   $\Box$

\bigskip
{\bf Proof of Theorem \ref{thDeviation}.} Let $\delta>0$ be a small number. Denote $T=\delta^{-2}$,
$m=m_n=[Tnx_n^2]$, $l=l_n=[n/m]$. Then $l\sim \frac{1}{Tx_n^2}$, $n x_n^2/m\sim T^{-1}$. Let $\phi(y)$ be a non-decreasing  Lipschitz function such that $I\{y\le 1\}\le \phi(y)\le I\{y\le 1+\delta/32\}$. Then by the definition of independence,
\begin{align*}
 \Capc\Big( &\max_{k\le n}|S_k|\le \sqrt{n}x_n\Big)\le \Capc\Big(\max_{j\le l}\max_{m(j-1)<k\le mj}|S_k|\le \sqrt{n}x_n\Big)\\
\le & \Sbep\left[\prod_{j=1}^l\phi\left(\max_{m(j-1)<k\le mj}|S_k|/(\sqrt{n}x_n)\right)\right]\\
= & \Sbep\left[\prod_{j=1}^{l-1}\phi\left(\max_{m(j-1)<k\le mj}|S_k|/(\sqrt{n}x_n)\right)\right.\\
& \qquad\quad \left. \cdot\Sbep\left[\phi\left(\max_{m(l-1)<k\le ml}|S_k-S_{m(l-1)}+S_{m(l-1)}|/(\sqrt{n}x_n)\right)\Big|S_k:k\le m(l-1)\right]\right].
\end{align*}
When $\phi\left(\max_{m(j-2)<k\le m(j-1)}|S_k|/(\sqrt{n}x_n)\right)\ne 0$, we have  $|S_{m(j-1)}|\le (1+\delta/32)\sqrt{n}x_n$, and then
\begin{align*}
&\Sbep\left[\phi\left(\max_{m(j-1)<k\le mj}|S_k-S_{m(j-1)}+S_{m(j-1)}|/(\sqrt{n}x_n)\right)\Big|S_k:k\le m(j-1)\right]\\
&\quad \le    \sup_{|y|\le (1+\delta/32)\sqrt{n}x_n}
\Sbep\left[\phi\left(\max_{ k\le m}|S_k+y|/(\sqrt{n}x_n)\right) \right]\\
&\quad \le  \sup_{|y|\le 2 T^{-1/2}}
\Sbep\left[\phi\left(\frac{\sup_{0\le t\le 1}|W_m(t)+y|}{\sqrt{n}x_n/\sqrt{m}}\right) \right],
\end{align*}
 by noting that $\{S_k-S_{m(j-1)}; k=m(j-1)+1,\ldots, mj\}$ and $\{S_k; k=1,\ldots, m\}$ are   identically distributed under $\Sbep$.
Note $\sqrt{n}x_n /\sqrt{m}\to T^{-1/2}$. Choose $h_m(x,y)= \phi\left(\frac{\sup_{0\le t\le 1}|x(t)+y(t)|}{\sqrt{n}x_n/\sqrt{m}}\right)$, $h(x,y)= \phi\left(\sup_{0\le t\le 1}|x(t)+y(t)|/T^{-1/2}\right)$
and $K=\{y(t)\equiv y: |y|\le 2 T^{-1/2}\}$ in (\ref{eqthFCLT.2}). By Corollary \ref{cor1}, uniformly in $|y|\le 2 T^{-1/2}$,
\begin{align*}
 \lim_{n\to \infty} \Sbep &\left[\phi\left(\frac{\sup_{0\le t\le 1}|W_m(t)+y|}{\sqrt{n}x_n/\sqrt{m}}\right) \right] \\
=& \Sbep\left[\phi\left(\sup_{0\le t\le 1}|W(t)+y|/T^{-1/2}\right) \right] \\
\le &  \Capc\Big(\sup_{0\le t\le 1}|W(t)+y|\le (1+\delta/32)T^{-1/2}\Big) \\
\le & \frac{4}{\pi}\exp\Big\{-\frac{\pi^2 T\overline{\sigma}^2}{8x^2(1+\delta/32)^2}\Big\}\;\; (\text{ by Lemma \ref{lemma7}}).
\end{align*}
So, there is a $n_0$ such that for all $n\ge n_0$,
\begin{align*}
  \sup_{|y|\le 2 T^{-1/2}}\Sbep &\left[\phi\left(\frac{\sup_{0\le t\le 1}|W_m(t)+y|}{\sqrt{n}x_n/\sqrt{m}}\right) \right] \\
 \le & \frac{8}{\pi}\exp\Big\{-\frac{\pi^2 T\overline{\sigma}^2}{8x^2(1+\delta/32)^2}\Big\} \le
 \exp\Big\{-\frac{\pi^2 T\overline{\sigma}^2}{8x^2(1+\delta/8)^2}\Big\}
\end{align*}
if $T=\delta^{-2}$ is large enough.  Hence
$$ \log \Capc\Big(\max_{k\le n}|S_k|\le \sqrt{n}x_n\Big)\le  -l \frac{\pi^2 T\overline{\sigma}^2}{8(1+\delta/8)^2}\sim -\frac{\pi^2   \overline{\sigma}^2}{8(1+\delta/8)^2}x_n^{-2}. $$
 It follows that
 $$ \limsup_{n\to \infty} x_n^2\log \Capc\Big(\max_{k\le n}|S_k|\le \sqrt{n}x_n\Big)
 \le -\frac{\pi^2  \overline{  \sigma}^2}{8}. $$

 Next, we consider the lower bound. Let   $\phi$ be defined as above, and $\phi_1$ be a non-decreasing  Lipschitz function such that $I\{y\le 1-\delta/32\}\le \phi(y)\le I\{y\le 1\}$. Let $\epsilon>0$ be a number whose value will be given later. Then
 \begin{align*}
&\Capc\Big(\max_{k\le n}|S_k|\le \sqrt{n}x_n\Big)\ge \Capc\Big(\max_{j\le l+1}\max_{m(j-1)<k\le mj}|S_k|\le \sqrt{n}x_n\Big)\\
\ge & \Sbep\left[\prod_{j=1}^{l+1}\phi_1\left(\max_{m(j-1)<k\le mj}|S_k|/(\sqrt{n}x_n)\right)\phi\left(|S_{jm}|/(\epsilon\sqrt{n}x_n)\right)\right]\\
= & \Sbep\left[\prod_{j=1}^{l}\phi_1\left(\frac{\max_{m(j-1)<k\le mj}|S_k|}{\sqrt{n}x_n}\right)\phi\left(\frac{|S_{jm}|}{\epsilon\sqrt{n}x_n}\right)\right.\\
&  \quad \left. \cdot\Sbep\left[\phi_1\left(\frac{\max_{ml<k\le m(l+1)}|S_k|}{\sqrt{n}x_n}\right)\phi\left(\frac{|S_{(l+1)m}|}{\epsilon\sqrt{n}x_n}\right)\Big|S_k:k\le m(l-1)\right]\right].
\end{align*}
When $\phi\left(\frac{|S_{(j-1)m}|}{\epsilon\sqrt{n}x_n}\right)\ne 0$, we have $|S_{(j-1)m}|\le  (1+\delta/32)\epsilon \sqrt{n}x_n$,  and then
\begin{align*}
 \Sbep &\left[\phi_1\left(\frac{\max_{m(j-1)<k\le mj}|S_k|}{\sqrt{n}x_n}\right)\phi
\left(\frac{|S_{jm}|}{\epsilon\sqrt{n}x_n}\right)\Big|S_k:k\le m(j-1)\right]\\
\ge  & \inf_{|y|\le (1+\delta/32)\epsilon \sqrt{n}x_n}\Sbep\left[\phi_1\left(\frac{\max_{k\le m}|S_k+y|}{\sqrt{n}x_n}\right)\phi
\left(\frac{|S_m+y|}{\epsilon\sqrt{n}x_n}\right) \right]\\
\ge & \inf_{|y|\le 2\epsilon T^{-1/2}}\Sbep\left[\phi_1\left(\frac{\sup_{0\le t\le 1} |W_n(t)+y|}{\sqrt{n}x_n/\sqrt{m}}\right)\phi
\left(\frac{|W_m(1)+y|}{\epsilon\sqrt{n}x_n/\sqrt{m}}\right) \right].
\end{align*}
By Corollary \ref{cor1}, uniformly in $|y|\le 2\epsilon T^{-1/2}$,
\begin{align*}
\lim_{n\to\infty}& \Sbep  \left[\phi_1\left(\frac{\sup_{0\le t\le 1} |W_n(t)+y|}{\sqrt{n}x_n/\sqrt{m}}\right)\phi
\left(\frac{|W_m(1)+y|}{\epsilon\sqrt{n}x_n/\sqrt{m}}\right) \right]\\
=& \widetilde{\mathbb E}  \left[\phi_1\left(\frac{\sup_{0\le t\le 1} |W(t)+y|}{T^{-1/2}}\right)\phi
\left(\frac{|W(1)+y|}{\epsilon T^{-1/2}}\right) \right] \\
\ge & \widetilde{\Capc}\Big(\sup_{0\le t\le 1} |W(t)+y|\le (1-\delta/32)T^{-1/2}, |W(1)+y|\le  \epsilon T^{-1/2}\Big).
\end{align*}
 Choose $\epsilon =\frac{\delta}{128}$. By Lemma \ref{lemma9}, if $T=\delta^{-2}$ is large enough,
\begin{align*}
  T^{-1}\log & \left(\inf_{|y|\le 2\epsilon T^{-1/2}}\widetilde{\Capc}\Big(\sup_{0\le t\le 1} |W(t)+y|\le (1-\delta/32)T^{-1/2}, |W(1)+y|\le  \epsilon T^{-1/2}\Big)\right)\\
 & \ge  -\frac{T\pi^2\underline{\sigma}^2}{8}( 1-\delta/32-4\epsilon)^{-2} (1-\delta/16)^{-1}
 \ge -\frac{T\pi^2\underline{\sigma}^2}{8}  (1-\delta/16)^{-3}.
 \end{align*}
 Hence, there is a $n_0$ such that for all $n\ge n_0$,
 \begin{align*}
 \Sbep   \left[\phi\left(\frac{\max_{m{j-1}<k\le mj}|S_k|}{\sqrt{n}x_n}\right)\phi_1
\left(\frac{|S_{jm}|}{\epsilon\sqrt{n}x_n}\right)\Big|S_k:k\le m(j-1)\right]
 \ge \exp \left\{-  \frac{T\pi^2\underline{\sigma}^2}{8}( 1-\delta/16)^{-4}\right\}.
\end{align*}
It follows that
$$\log \Capc\Big(\max_{k\le n}|S_k|\le \sqrt{n}x_n\Big)\ge -(l+1)T \frac{\pi^2\underline{\sigma}^2}{8}( 1-\delta/16)^{-4}\sim -\frac{\pi^2\underline{\sigma}^2}{8}( 1-\delta/16)^{-4}x_n^{-2}. $$
 Hence,
 $$ \liminf_{n\to \infty} x_n^2\log \Capc\Big(\max_{k\le n}|S_k|\le \sqrt{n}x_n\Big)\ge
 -\frac{\pi^2\underline{\sigma}^2}{8}. $$
 The proof of (\ref{eqDeviation.1}) is completed. The proof of (\ref{eqDeviation.1}) is similar. $\Box$.


\end{document}